\documentclass[aos,seceqn,citesort,dvips]{arximspdf}
\usepackage{dcolumn}

\doi{10.1214/09-AOS747}
\volume{38}
\issue{6}
\pubyear{2010}
\firstpage{3387}
\lastpage{3411}

\makeatletter

\newcolumntype{d}[1]{D{.}{.}{#1}}

\newtheorem{theorem}{Theorem}[section]
\newtheorem{lemma}{Lemma}
\newtheorem{proposition}{Proposition}

\newproclaim{assumption}{Assumption}

\makeatother

\begin{document}
\begin{frontmatter}

\title{Sequentially interacting Markov chain Monte Carlo methods}
\runtitle{Sequentially interacting MCMC}

\begin{aug}
\author[A]{\fnms{Anthony} \snm{Brockwell}\ead[label=e1]{Abrock@stat.cmu.edu}},
\author[B]{\fnms{Pierre} \snm{Del Moral}\ead[label=e2]{Pierre.Del-Moral@inria.fr}} and
\author[C]{\fnms{Arnaud} \snm{Doucet}\corref{}\ead[label=e3]{Arnaud@stat.ubc.ca}}
\runauthor{A. Brockwell, P. Del Moral and A. Doucet}
\affiliation{Carnegie Mellon University, INRIA Bordeaux Sud-Ouest
and~University~of~British Columbia}
\address[A]{A. Brockwell\\
Department of Statistics\\
Carnegie Mellon University\\
5000 Forbes Avenue\\
Pittsburgh, Pennsylvania 15123\\
USA\\
\printead{e1}}
\address[B]{P. Del Moral\\
Centre INRIA Bordeaux Sud-Ouest\\
Institut de Math\'{e}matiques\\
Universit\'{e} Bordeaux I\\
33405 Talence Cedex\\
France\\
\printead{e2}}
\address[C]{A. Doucet\\
Departments of Computer Science\\
\quad and Statistics\\
University of British Columbia\\
333-6356 Agricultural Road\\
Vancouver, BC, V6T 1Z2\\
Canada\\
\printead{e3}}
\end{aug}

% HISTORY:
\received{\smonth{1} \syear{2008}}
\revised{\smonth{8} \syear{2009}}

% ABSTRACT
%
\begin{abstract}
Sequential Monte Carlo (SMC) is a methodology for sampling
approximately from a sequence of probability distributions of
increasing dimension and estimating their normalizing constants. We
propose here an alternative methodology named Sequentially Interacting
Markov Chain Monte Carlo (SIMCMC). SIMCMC methods work by generating
interacting non-Markovian sequences which behave asymptotically like
independent Metropolis--Hastings (MH) Markov chains with the desired
limiting distributions. Contrary to SMC, SIMCMC allows us to
iteratively improve our estimates in an MCMC-like fashion. We establish
convergence results under realistic verifiable assumptions and
demonstrate its performance on several examples arising in Bayesian
time series analysis.
\end{abstract}

% KEYWORDS
%
\begin{keyword}[class=AMS]
\kwd[Primary ]{65C05}
\kwd{60J05}
\kwd[; secondary ]{62F15}.
\end{keyword}
\begin{keyword}
\kwd{Markov chain Monte Carlo}
\kwd{normalizing constants}
\kwd{sequential Monte Carlo}
\kwd{state-space models}.
\end{keyword}

\end{frontmatter}

%s1 ###
\section{Introduction}\label{sec1}

Let us consider a sequence of probability distributions $\{\pi_{n}\}
_{n\in
\mathbb{T}}$ where $\mathbb{T}= \{ 1,2,\ldots,P \} $, which we will
refer to as ``target'' distributions. We
shall also refer to $n$ as the time index. For ease of presentation, we
shall assume here that $\pi_{n} ( d\mathbf{x}_{n} ) $ is defined
on a measurable space $ ( E_{n},\mathcal{F}_{n} ) $ where $E_{1}=E$,
$\mathcal{F}_{1}=\mathcal{F}$ and $E_{n}=E_{n-1}\times E$, $\mathcal
{F}_{n}=%
\mathcal{F}_{n-1}\times\mathcal{F}$. We denote $\mathbf{x}_{n}= (
x_{1},\ldots,x_{n} ) $ where $x_{i}\in E$ for $i=1,\ldots,n$. Each $\pi
_{n} ( d\mathbf{x}_{n} ) $ is assumed to admit a density $\pi
_{n} ( \mathbf{x}_{n} ) $ with respect to a $\sigma$-finite
dominating measure denoted $d\mathbf{x}_{n}$ and $d\mathbf
{x}_{n}=d\mathbf{x}%
_{n-1}\times dx_{n}$. Additionally, we have
\[
\pi_{n} ( \mathbf{x}_{n} ) =\frac{\gamma_{n} ( \mathbf{x}%
_{n} ) }{Z_{n}},
\]
where $\gamma_{n}\dvtx E_{n}\rightarrow\mathbb{R}^{+}$ is known
pointwise and
the normalizing constant $Z_{n}$ is unknown.

In a number of important applications, it is desirable to be able to sample
from the sequence of distributions $\{\pi_{n}\}_{n\in\mathbb{T}}$ and to
estimate their normalizing constants $\{Z_{n}\}_{n\in\mathbb{T}}$; the most
popular statistical application is the class of nonlinear non-Gaussian
state-space models detailed in Section~\ref{sec:applications}. In this
context, $\pi_{n}$ is the posterior distribution of the hidden state
variables from time $1$ to $n$ given the observations from time $1$ to $n$
and $Z_{n}$ is the marginal likelihood of these observations. Many other
applications, including contingency tables and population genetics, are
discussed in~\cite{Delm04,Douc01} and~\cite{Liu2001}.

A standard approach to solve this class of problems relies on Sequential
Monte Carlo (SMC) methods; see~\cite{Douc01} and~\cite{Liu2001} for a review
of the literature. In the SMC approach, the target distributions are
approximated by a large number of random samples, termed particles, which
are carried forward over time by using a combination of sequential
importance sampling and resampling steps. These methods have become the
tools of choice for sequential Bayesian inference but, even when there
is no
requirement for ``real-time'' inference, SMC algorithms are
increasingly used
as an alternative to MCMC; see, for example, \cite
{Chopin2002,DelMoralDoucetJasra2006}
and~\cite{Liu2001} for applications to
econometrics models, finite mixture models and contingency tables. They also
allow us to implement goodness-of-fit tests easily in a time series context
whereas a standard MCMC implementation is cumbersome~\cite{gerlach1999}.
Moreover, they provide an estimate of the marginal likelihood of the data.

The SMC methodology is now well established and many theoretical convergence
results are available~\cite{Delm04}. Nevertheless, in practice, it is
typically impossible to, a priori, determine the number of particles
necessary to achieve a fixed precision for a given application. In such
scenarios, users typically perform multiple runs for an increasing
number of
particles until stabilization of the Monte Carlo estimates is observed.
Moreover, SMC algorithms are substantially different from MCMC algorithms
and can appear difficult to implement for nonspecialists.

In this paper, we propose an alternative to SMC named \textit{Sequentially
Interacting Markov Chain Monte Carlo} (SIMCMC). SIMCMC methods allow us to
compute Monte Carlo estimates of the quantities of interest iteratively as
they are, for instance, when using MCMC methods. This allows us to refine
the Monte Carlo estimates until a suitably chosen stopping time.
Furthermore, for people familiar with MCMC methods, SIMCMC methods are
somewhat simpler than SMC methods to implement, because they only rely
on MH
steps. However, SIMCMC methods are not a class of MCMC methods. These are
non-Markovian algorithms which can be interpreted as an approximation
of $P$
``ideal'' standard MCMC chains. It is based on the same key idea as SMC
methods; that is as $\pi_{n+1} ( \mathbf{x}_{n} ) =\int\pi
_{n+1} ( \mathbf{x}_{n+1} ) \,dx_{n+1}$ is often ``close'' to $\pi
_{n} ( \mathbf{x}_{n} ) $, it is sensible to use $\pi_{n} (
\mathbf{x}_{n} ) $ as part of a proposal distribution to sample $\pi
_{n+1} ( \mathbf{x}_{n+1} ) $. In SMC methods, the correction
between the proposal distribution and the target distribution is performed
using Importance Sampling whereas in SIMCMC methods it is performed using
an MH step. Such a strategy is computationally much more efficient than
sampling separately from each target distribution using standard MCMC
methods and also provides direct estimates of the normalizing constants
$%
\{Z_{n}\}_{n\in\mathbb{T}}$.

The potential real-time applications of SIMCMC methods are also worth
commenting on. SMC methods have been used in various real-time engineering
applications, for example, in neural decoding~\cite{brockwell2004} and in
target tracking~\cite{Gilks99,septier2009}. In these problems,
it is
important to be able to compute functionals of the posterior distributions
of some quantity of interest, but it must also be done in real-time. SMC
methods work with collections of particles that are updated
sequentially to
reflect these distributions. Clearly, in such real-time problems it is
important that the collections of particles are not too large, or else the
computational burden can cause the SMC algorithm to fall behind the system
being analyzed. SIMCMC methods provide a very convenient way to make optimal
use of what computing power is available. Since SIMCMC works by adding one
particle at a time to collections representing distributions, we can simply
run it continually in between arrival of successive observations, and it
will accrue as many particles as it can in whatever amount of time is taken.

Related ideas where we also have a sequence of nested MCMC-like chains
``feeding'' each other and targeting a sequence of increasingly complex
distributions have recently appeared in statistics~\cite{kou2006} and
physics~\cite{lyman2006}. In the equi-energy sampler~\cite{kou2006}, the
authors consider a sequence of distributions indexed by a temperature
and an
energy truncation whereas in~\cite{lyman2006} the authors consider a
sequence of coarse-grained distributions. It is also possible to think of
SIMCMC methods and the algorithms in~\cite{kou2006} and \cite
{lyman2006} as
nonstandard adaptive MCMC schemes \cite
{andrieumoulines2006,robertsrosenthal2007}
where the parameters to be adapted are probability
distributions instead of finite-dimensional parameters. Our convergence
results rely partly on ideas developed in this field~\cite{andrieumoulines2006}.

The rest of the paper is organized as follows. In Section~\ref{sec:SIMCMC},
we describe \mbox{SIMCMC} methods, give some guidelines for the design of
efficient algorithms and discuss implementation issues. In Section \ref%
{sec:convergence}, we present some convergence results. In Section \ref%
{sec:applications}, we demonstrate the performance of this algorithm for
various Bayesian time series problems and compare it to SMC. Finally, we
discuss a number of further potential extensions in Section \ref%
{sec:discussion}. The proofs of the results in Section~\ref{sec:convergence}
can be found in \hyperref[app]{Appendix}.

%s2 ###
\section{Sequentially interacting Markov chain Monte Carlo}\label{sec:SIMCMC}

%s2.1 ###
\subsection{The SIMCMC algorithm}

Let $i$ be the iteration counter. The SIMCMC algorithm constructs $P$
sequences $ \{ \mathbf{X}{_{1}^{ ( i ) }} \} , \{
\mathbf{X}{_{2}^{ ( i ) }} \} ,\ldots, \{ \mathbf{X}{%
_{P}^{ ( i ) }} \} $. In Section~\ref{sec:convergence}, we
establish weak necessary conditions ensuring that as $i$ approaches
infinity, the distribution of $\mathbf{X}{_{n}^{ ( i ) }}$
approaches $\pi_{n}$; we will assume here that these conditions are
satisfied to explain the rationale behind our algorithm. To specify the
algorithm, we require a sequence of $P$ proposal distributions,
specified by
their densities
\[
q_{1}(x_{1}),q_{2}(\mathbf{x}_{1},x_{2}),\ldots,q_{P}(\mathbf{x}%
_{P-1},x_{P}).
\]
Each $q_{n}\dvtx E_{n-1}\times E\rightarrow\mathbb{R}^{+}$
($E_{-1}=\varnothing$%
) is a probability density in its last argument $x_{n}$ with respect to
$%
dx_{n}$, which may depend (for $n=2,\ldots,P$) on the first argument.
Proposals are drawn from $q_{1}(\cdot)$ for updates of the sequence $\{
\mathbf{X}{_{1}^{ ( i ) }}\}$, from $q_{2}(\cdot)$ for updates of
the sequence $\{\mathbf{X}{_{2}^{ ( i ) }}\}$, and so on.
(Selection of proposal distributions is discussed in Section \ref%
{sec:settings}.) Based on these proposals, we define the weights
%
%e2.1 ###
%
\begin{eqnarray}\label{eq:unnormalizedweights}
{w}_{1} ( \mathbf{x}_{1} ) & = & \frac{\gamma_{1} ( \mathbf{x}%
_{1} ) }{q_{1} ( \mathbf{x}_{1} ) }, \nonumber\\[-8pt]\\[-8pt]
{w}_{n} ( \mathbf{x}_{n} ) & = & \frac{\gamma_{n} ( \mathbf{x}%
_{n} ) }{\gamma_{n-1} ( \mathbf{x}_{n-1} ) q_{n} (
\mathbf{x}_{n-1},x_{n} ) },\qquad n=2,\ldots,P.\nonumber
\end{eqnarray}
For any measure $\mu_{n-1}$ on $ ( E_{n-1},\mathcal{F}_{n-1} ) $,
we define
\[
( \mu_{n-1}\times q_{n} ) ( d\mathbf{x}_{n} ) =\mu
_{n-1} ( d\mathbf{x}_{n-1} ) q_{n} ( \mathbf{x}%
_{n-1},dx_{n} )
\]
and
%
%e2.2 ###
%
\begin{equation} \label{eq:supporttargets}
\mathcal{S}_{n}= \{ \mathbf{x}_{n}\in E_{n}\dvtx\pi_{n} ( \mathbf{x}%
_{n} ) >0 \} .
\end{equation}
We also denote by $\widehat{\pi}_{n}^{ ( i ) }$ the empirical
measure approximation of the target distribution $\pi_{n}$ given by
%
%e2.3 ###
%
\begin{equation} \label{eq:empiricaldistributiontarget}
\widehat{\pi}_{n}^{ ( i ) } ( d\mathbf{x}_{n} ) =\frac{1}{%
i+1}\sum_{m=0}^{i}\delta_{\mathbf{X}{_{n}^{ ( m ) }}} ( d%
\mathbf{x}_{n} ) .
\end{equation}

Intuitively, the SIMCMC algorithm proceeds as follows. At each
iteration $i$
of the algorithm, the algorithm samples $\mathbf{X}{_{n}^{ ( i ) }}$
for $n\in\mathbb{T}$ by first sampling $\mathbf{X}{_{1}^{ ( i ) }}$%
, then $\mathbf{X}{_{2}^{ ( i ) }}$ and so on. For $n=1$, $ \{
\mathbf{X}{_{1}^{ ( i ) }} \} $ is a standard Markov chain
generated using an independent MH sampler of invariant distribution $\pi
_{1} ( \mathbf{x}{_{1}} ) $ and proposal distribution $q_{1} (
\mathbf{x}{_{1}} ) $. For $n=2$, we would like to approximate an
independent MH sampler of invariant distribution $\pi_{2} ( \mathbf{x}{
_{2}} ) $ and proposal distribution $ ( \pi_{1}\mathbf{\times}%
q_{2} ) ( \mathbf{x}{_{2}} ) $. As it is impossible to sample
from $\pi_{1}$ exactly, we replace $\pi_{1}$ at iteration $i$ by its
current empirical measure approximation $\widehat{\pi}_{1}^{ ( i )
}$. Similarly for $n>2$, we approximate an MH sampler of invariant
distribution $\pi_{n} ( \mathbf{x}{_{n}} ) $ and proposal
distribution $ ( \pi_{n-1}\mathbf{\times}q_{n} ) ( \mathbf{x}%
_{n} ) $ by replacing $\pi_{n-1}$ at iteration $i$ by its current
empirical measure approximation $\widehat{\pi}_{n-1}^{ ( i ) }$.
The sequences $ \{ \mathbf{X}{_{2}^{ ( i ) }} \} ,\ldots
, \{ \mathbf{X}{_{P}^{ ( i ) }} \} $ generated this way
are clearly non-Markovian.
\newpage

\noindent\hrulefill%

\begin{center}
\textbf{Sequentially interacting Markov chain Monte Carlo}
\end{center}

$\bullet$ \textsf{Initialization, }$i=0$

\qquad$\bullet$ \textsf{For} $n\in\mathbb{T}$\textsf{, set randomly }$%
\mathbf{X}_{n}^{(0)}=\mathbf{x}{_{n}^{ ( 0 ) }\in}$ $\mathcal{S}%
_{n}$\textsf{.}

$\bullet$ \textsf{For iteration }$i\geq1$

\qquad$ \bullet$ \textsf{For }$n=1$

\qquad\qquad$\bullet$ \textsf{Sample }$\mathbf{X}_{1}^{\ast(
i ) }\sim q_{1} ( \cdot) $\textsf{.}

\qquad\qquad$\bullet$ \textsf{With probability }
%
%e2.4 ###
%
\begin{equation} \label{eq:MHratio1}
\alpha_{1}\bigl(\mathbf{X}_{1}^{(i-1)},\mathbf{X}_{1}^{\ast( i )
}\bigr)=1\wedge\frac{w_{1} ( \mathbf{X}_{1}^{\ast( i ) } ) }{%
w_{1} ( \mathbf{X}_{1}^{(i-1)} ) }
\end{equation}
\textsf{ }

\qquad\qquad\textsf{set }$\mathbf{X}_{1}^{(i)}=\mathbf{X}_{1}^{{\ast} (
i ) }$\textsf{, otherwise set }$\mathbf{X}_{1}^{(i)}=\mathbf{X}%
_{1}^{(i-1)}$.

\qquad$ \bullet$ \textsf{For }$n=2,\ldots,P$

\qquad\qquad$\bullet$ \textsf{Sample }$\mathbf{X}_{n}^{\ast(
i ) }\sim( \widehat{\pi}_{n-1}^{(i)}\times q_{n} ) (
\cdot) $\textsf{.}

\qquad\qquad$\bullet$ \textsf{With probability}
%
%e2.5 ###
%
\begin{equation} \label{eq:MHratiomore}
\alpha_{n}\bigl(\mathbf{X}_{n}^{(i-1)},\mathbf{X}_{n}^{\ast( i )
}\bigr)=1\wedge\frac{w_{n} ( \mathbf{X}_{n}^{\ast( i ) } ) }{%
w_{n} ( \mathbf{X}_{n}^{(i-1)} ) }
\end{equation}

\qquad\qquad\textsf{set} $\mathbf{X}_{n}^{(i)}=\mathbf{X}_{n}^{\ast(
i ) }$\textsf{, otherwise set} $\mathbf{X}_{n}^{(i)}=\mathbf{X}%
_{n}^{(i-1)}$.

\noindent\hrulefill\vspace*{6pt}%

The (ratio of) normalizing constants can easily be estimated by
%
%e2.6 ###
%
\begin{eqnarray}\label{eq:estimatesofnormalizingconstant}
\widehat{Z}{}^{ ( i ) }_{1}& = &\frac{1}{i}\sum_{m=1}^{i}w_{1} \bigl(
\mathbf{X}_{1}^{\ast( m ) } \bigr) , \nonumber\\[-8pt]\\[-8pt]
\widehat{ \biggl( \frac{Z_{n}}{Z_{n-1}} \biggr) }^{ ( i ) }& = &\frac
{1}{%
i}\sum_{m=1}^{i}w_{n} \bigl( \mathbf{X}_{n}^{\ast( m ) } \bigr)
.\nonumber
\end{eqnarray}
Equation (\ref{eq:estimatesofnormalizingconstant}) follows from the identity
\[
\frac{Z_{n}}{Z_{n-1}}=\int{w}_{n} ( \mathbf{x}_{n} ) ( \pi
_{n-1}\times q_{n} ) ( d\mathbf{x}_{n} )
\]
and the fact that asymptotically (as $i\rightarrow\infty$) $\mathbf{X}
_{n}^{\ast( i ) }$ is distributed according to $ ( \pi_{n-1}%
\mathbf{\times}q_{n} ) ( \mathbf{x}_{n} ) $ under weak
conditions given in Section~\ref{sec:convergence}.

%s2.2 ###
\subsection{Algorithm settings}\label{sec:settings}

In a similar manner to SMC methods, the performance of the SIMCMC algorithm
depends on the proposal distributions. However, it is possible to devise
some useful guidelines for this sequence of \mbox{(pseudo-)}independent
samplers,
using reasoning similar to that adopted in the SMC framework.
Asymptotically, $\mathbf{X}_{n}^{\ast( i ) }$ is distributed
according to $ ( \pi_{n-1}\mathbf{\times}q_{n} ) ( \mathbf{x}%
_{n} ) $ and $w_{n}(\mathbf{x}_{n})$ is just the importance weight (up
to a normalizing\vadjust{\goodbreak} constant) between $\pi_{n}(\mathbf{x}_{n})$ and $ (
\pi_{n-1}\mathbf{\times}q_{n} ) ( \mathbf{x}_{n} ) $. The
proposal distribution minimizing the variance of this importance weight is
simply given by
%
%e2.7 ###
%
\begin{equation} \label{eq:optimalproposal}
q_{n}^{\mathrm{opt}} ( \mathbf{x}_{n-1},x_{n} ) =\overline{\pi}%
_{n} ( \mathbf{x}_{n-1},x_{n} ),
\end{equation}
where $\overline{\pi}_{n} ( \mathbf{x}_{n-1},x_{n} ) $ is the
conditional density of $x_{n}$ given $\mathbf{x}_{n-1}$ under $\pi_{n}$,
that is,
%
%e2.8 ###
%
\begin{equation}\label{eq:conditionaldensitypin}
\overline{\pi}_{n} ( \mathbf{x}_{n-1},x_{n} ) =\frac{\pi
_{n} ( \mathbf{x}_{n} ) }{\pi_{n} ( \mathbf{x}_{n-1} ) }.
\end{equation}
In the SMC literature, $\overline{\pi}_{n} ( \mathbf{x}%
_{n-1},x_{n} ) $ is called the ``optimal'' importance
density~\cite{Douc00}. This yields
%
%e2.9 ###
%
\begin{equation} \label{eq:optimalweight}
w_{n}^{\mathrm{opt}} ( \mathbf{x}_{n} ) \propto\pi_{n/n-1} (
\mathbf{x}_{n-1} ),
\end{equation}
where
%
%e2.10 ###
%
\begin{equation}\label{eq:optimalimportanceweight}
\pi_{n/n-1} ( \mathbf{x}_{n-1} ) =\frac{\pi_{n} ( \mathbf{x}%
_{n-1} ) }{\pi_{n-1} ( \mathbf{x}_{n-1} ) }
\end{equation}
with
\[
\pi_{n} ( \mathbf{x}_{n-1} ) =\int_{E}\pi_{n} ( \mathbf{x}%
_{n} ) \,dx_{n}.
\]
In this case, as $w_{n}^{\mathrm{opt}} ( \mathbf{x}_{n} ) $ is
independent of $x_{n}$, the algorithm described above can be further
simplified. It is indeed possible to decide whether to accept or reject a
candidate even before sampling it. This is more computationally efficient
because if the move is to be rejected there is no need to sample the
candidate. In most applications, it will be difficult to sample from
(\ref%
{eq:optimalproposal}) and/or to compute (\ref{eq:optimalweight}) as it
involves computing $\pi_{n} ( \mathbf{x}_{n-1} ) $ up to a
normalizing constant. In this case, we recommend approximating (\ref%
{eq:optimalproposal}). Similar strategies have been developed successfully
in the SMC framework~\cite{fearnhead1999,Douc00,liuchen1998,Pitt99}.
The advantages of such sampling strategies in the SIMCMC
case will be demonstrated in the simulation section.

Generally speaking, most of the methodology developed in the SMC setting
can be directly reapplied here. This includes the use of
Rao-Blackwellisa\-tion techniques to reduce the dimensionality of the target
distributions~\cite{Douc00,liuchen1998} or of auxiliary
particle-type ideas where we build target distributions biased toward
``promising'' regions of the space~\cite{fearnhead1999,Pitt99}.

%s2.3 ###
\subsection{Implementation issues}

%s2.3.1 ###
\subsubsection{Burn-in and storage requirements}
\label{subsubsec:burnin}

We have presented the algorithm without any burn-in. This can be easily
included if necessary by considering at iteration $i$ of the algorithm
\[
\widehat{\pi}_{n}^{(i)}(d\mathbf{x}_{n})={\frac{1}{i+1-l(i,B)}}%
\sum_{m=l(i,B)}^{i}\delta_{\mathbf{X}_{n}^{ ( m ) }}(d\mathbf{x}%
_{n}),
\]
where
\[
l(i,B)=0\vee\bigl( ( i-B ) \wedge B\bigr),\vadjust{\goodbreak}
\]
where $B$ is an appropriate number of initial samples to be discarded as
burn-in. Note that when $i\geq2B$, we have $l(i,B)=B$.

Note that in its original form, the SIMCMC algorithm requires storing the
sequences $ \{ \mathbf{X}{_{n}^{ ( i ) }} \} _{n\in
\mathbb{T}}$. This could be expensive if the number of target
distributions $%
P$ and/or the number of iterations of the SIMCMC are large. However, in many
scenarios of interest, including nonlinear non-Gaussian state-space models
or the scenarios considered in~\cite{DelMoralDoucetJasra2006}, it is
possible to drastically reduce these storage requirements as we are only
interested in estimating the marginals $ \{ \pi_{n} ( x{_{n}}%
) \} $ and we have $w_{n} ( \mathbf{x}_{n} )
=w_{n} ( x_{n-1},x_{n} ) $ and $q_{n} ( \mathbf{x}%
_{n-1},x_{n} ) =q_{n} ( x_{n-1},x_{n} ) $. In such cases, we
only need to store $ \{ X{_{n}^{ ( i ) }} \} _{n\in
\mathbb{T}}$, resulting in significant memory savings.

%s2.3.2 ###
\subsubsection{Combining sampling strategies}

In practice, it is possible to combine the SIMCMC strategy with SMC
methods; that is we can generate say $N$ (approximate) samples from $%
\{ \pi_{n} \} _{n\in\mathbb{T}}$ using SMC then use the SIMCMC
strategy to increase the number of particles until the Monte Carlo estimates
stabilize. We emphasize that SIMCMC will be primarily useful in the context
where we do not have a predetermined computational budget. Indeed, if the
computational budget is fixed, then better estimates could be obtained by
switching the iteration $i$ and time $n$ loops in the SIMCMC
algorithm.\looseness=1

%s2.4 ###
\subsection{Discussion and extensions}\label{subsec:extensions}

Standard MCMC methods do not address the problem solved by SIMCMC methods.
Trans-dimensional MCMC methods~\cite{Green03} allow us to sample from a
sequence of ``related'' target distributions of different dimensions but
require the knowledge of the ratio of normalizing constants between these
target distributions. Simulated tempering and parallel tempering
require all
target distributions to be defined on the same space and rely on MCMC
kernels to explore each target distribution; see~\cite{jasra2006} for a
recent discussion of such techniques. As mentioned earlier in the
\hyperref[sec1]{Introduction}, ideas related to SIMCMC where a sequence
of ``ideal'' MCMC
algorithms is approximated have recently appeared in statistics \cite%
{kou2006} and physics~\cite{lyman2006}. However, contrary to these
algorithms, the target distributions considered here are of increasing
dimension and the proposed interacting mechanism is simpler. Whereas the
equi-energy sampler~\cite{kou2006} allows ``swap'' moves between
chains, we
only use the samples of the sequence $ \{ \mathbf{X}{_{n}^{ (
i ) }} \} $ to feed $ \{ \mathbf{X}{_{n+1}^{ ( i ) }}%
\} $ but $ \{ \mathbf{X}{_{n+1}^{ ( i ) }} \} $ is
never used to generate $ \{ \mathbf{X}{_{n}^{ ( i ) }} \} $.

There are many possible extensions of the SIMCMC algorithm. In this respect,
the SIMCMC algorithm is somehow a proof-of-concept algorithm demonstrating
that it is possible to make sequences targeting different distributions
interact without the need to define a target distribution on an extended
state space. For example, a simple modification of the SIMCMC algorithm
can be easily parallelized. Instead of sampling our candidate $\mathbf
{X}%
_{n}^{\ast( i ) }$ at iteration $i$ according to $(\widehat{\pi}%
_{n-1}^{(i)}\times q_{n}) ( \cdot) $ we can sample it
instead from $(\widehat{\pi}_{n-1}^{(i-1)}\times q_{n}) ( \cdot) $:
this allows us to simulate the sequences $ \{ \mathbf{X}{_{n}^{ (
i ) }} \} $ on $P$ parallel processors. It is straightforward to
adapt the convergence results given in Section~\ref{sec:convergence} to this
parallel version of SIMCMC.

In the context of real-time applications where $\pi_{n} ( \mathbf{x}%
_{n} ) $ is typically the posterior distribution $p (
\mathbf{x}_{n} \vert y_{1\dvtx n} ) $ of some states $\mathbf{x}_{n}$
given the observations $y_{1\dvtx n}$, \mbox{SIMCMC} methods can also
be very useful.
Indeed, SMC methods cannot easily address situations where the observations
arrive at random times whereas SIMCMC methods allow us to make optimal use
of what computing power is available by adding as many particles as possible
until the arrival of a new observation. In such cases, a standard
implementation would consist of updating our approximation of $\pi
_{n} ( \mathbf{x}_{n} ) $ at ``time'' $n$ by adding iteratively
particles to the approximations $\pi_{n-L+1} ( \mathbf{x}%
_{n-L+1} ) ,\ldots,\pi_{n-1} ( \mathbf{x}_{n-1} ) ,\pi
_{n} ( \mathbf{x}_{n} ) $ for a lag $L\geq1$ until the arrival of
data $y_{n+1}$. For $L=1$, such schemes have been recently proposed
independently in~\cite{septier2009}.

%s3 ###
\section{Convergence results}\label{sec:convergence}

We now present some convergence results for \mbox{SIMCMC}. Despite the
non-Markovian nature of this algorithm, we are here able to provide
realistic verifiable assumptions ensuring the asymptotic consistency of the
SIMCMC estimates. Our technique of proof rely on Poisson's equation
\cite{glynn1996}; similar tools have been used in \cite
{andrieumoulines2006} to
study the convergence of adaptive MCMC schemes and in
\cite{delmoralmiclo2003} to study the stability of self-interacting
Markov chains.

Let us introduce $B ( E_{n} ) =\{ f_{n}\dvtx E_{n}\rightarrow
\mathbb{R}$ such that $ \Vert f_{n} \Vert\leq1\} $
where $ \Vert f_{n} \Vert={\sup_{\mathbf{x}_{n}\in E_{n}}}%
\vert f_{n} ( \mathbf{x}_{n} ) \vert$. We denote\vspace*{-2pt} by $%
\mathbb{E}_{\mathbf{x}_{1\dvtx n}^{ ( 0 ) }} [ \cdot] $ the
expectation with respect to the distribution of the simulated
sequences initialized at $\mathbf{x}_{1\dvtx n}^{ ( 0 ) }:= (
\mathbf{x}_{1}^{ ( 0 ) },\mathbf{x}_{2}^{ ( 0 ) },\ldots,\break%
\mathbf{x}_{n}^{ ( 0 ) } ) $ and $\mathbb{N}_{0}=\mathbb{N}\cup
\{ 0 \} $. For any measure $\mu$ and test function $f$, we write
$\mu( f ) =\int\mu( dx ) f ( x ) $.

Our proofs rely on the following assumption.
\renewcommand{\theassumption}{A1}
\begin{assumption}\label{assumA1}
For any $n\in\mathbb{T}$, there exists $%
B_{n}<\infty$ such that for any $\mathbf{x}_{n}\in\mathcal{S}_{n}$%
%
%e3.1 ###
%
\begin{equation}\label{eq:incrementalweightbounded}
w_{n} ( \mathbf{x}_{n} ) \leq B_{n}.
\end{equation}
\end{assumption}

This assumption is quite weak and can be easily checked in all the examples
presented in Section~\ref{sec:applications}. Note that a similar assumption
also appears when $\mathbb{L}_{p}$ bounds are established for SMC methods
\cite{Delm04}.

Our first result establishes the convergence of the empirical averages
toward the correct expectations at the standard Monte Carlo rate.
\begin{theorem}
\label{theorem:convergenceofaverages} Given Assumption~\ref{assumA1},
for any $n\in\mathbb{T}$
and any $p\geq1$, there exist $C_{1,n},C_{2,p}<\infty$ such that for
any $%
\mathbf{x}_{1\dvtx n}^{ ( 0 ) }\in\mathcal{S}_{1}\times\cdots\times
\mathcal{S}_{n}$, $f_{n}\in B ( E_{n} ) $ and $i\in\mathbb{N}_{0}$%
\[
\mathbb{E}_{\mathbf{x}_{1\dvtx n}^{ ( 0 ) }} \bigl[ \bigl\vert\widehat{\pi}%
_{n}^{(i)} ( f_{n} ) -\pi_{n} ( f_{n} ) \bigr\vert^{p}%
\bigr] ^{1/p}\leq\frac{C_{1,n}C_{2,p}}{ ( i+1 ) ^{{1/2}}}.
\]
\end{theorem}

Using (\ref{eq:estimatesofnormalizingconstant}), a straightforward
corollary of Theorem~\ref{theorem:convergenceofaverages} is the following
result.
\begin{theorem}
\label{theorem:convergenceofnormalizingconstants}Given Assumption \ref
{assumA1}, for any
$n\in
\mathbb{T}$ and any $p\geq1$, there exist $C_{1,n},C_{2,p}<\infty$ such
that for any $\mathbf{x}_{1\dvtx n}^{ ( 0 ) }\in\mathcal{S}_{1}\times
\cdots\times\mathcal{S}_{n}$, $f_{n}\in B ( E_{n} ) $ and $i\in
\mathbb{N}_{0}$
\[
\mathbb{E}_{\mathbf{x}_{1}^{ ( 0 ) }} \bigl[ \bigl\vert\widehat{%
Z}{}^{ ( i ) }_{1}-Z_{1} \bigr\vert^{p} \bigr] ^{1/p}\leq\frac{%
B_{1}C_{1,1}C_{2,p}}{i^{{1/2}}},
\]
and for $n\in\mathbb{T}\setminus\{ 1 \} $
\[
\mathbb{E}_{\mathbf{x}_{1\dvtx n}^{ ( 0 ) }} \biggl[ \biggl\vert\widehat
{%
\biggl( \frac{Z_{n}}{Z_{n-1}} \biggr) }^{ ( i ) }-\frac{Z_{n}}{Z_{n-1}%
} \biggr\vert^{p} \biggr] ^{1/p}\leq\frac{B_{n}C_{1,n}C_{2,p}}{i^{{1/2}}}.
\]
\end{theorem}

By the Borel--Cantelli lemma, Theorems \ref
{theorem:convergenceofaverages} and
\ref{theorem:convergenceofnormalizingconstants} also ensure almost
sure convergence of the empirical averages and of the normalizing constant
estimates.

Our final result establishes that each sequence $ \{ \mathbf{X}%
_{n}^{ ( i ) } \} $ converges toward~$\pi_{n}$.

\begin{theorem}
\label{theorem:convergenceofmarginals}Given Assumption~\ref{assumA1},
for any $n\in\mathbb
{T}$, $%
\mathbf{x}_{1\dvtx n}^{ ( 0 ) }\in\mathcal{S}_{1}\times\cdots\times
\mathcal{S}_{n}$ and $f_{n}\in B ( E_{n} ) $ we have
\[
\lim_{i\rightarrow\infty}\mathbb{E}_{\mathbf{x}_{1\dvtx n}^{ (
0 ) }} \bigl[ f_{n} \bigl( \mathbf{X}_{n}^{ ( i ) } \bigr) %
\bigr] =\pi_{n} ( f_{n} ) .
\]
\end{theorem}

%s4 ###
\section{Applications}\label{sec:applications}

In this section, we will focus on the applications of \mbox
{SIMCMC} to nonlinear
non-Gaussian state-space models. Consider an unobserved $E$-valued Markov
process $ \{ X_{n} \} _{n\in\mathbb{T}}$ satisfying
\[
X_{1}\sim\mu( \cdot),\qquad X_{n} \vert\{
X_{n-1}=x \} \sim f ( x,\cdot) .
\]
We assume that we have access to observations $ \{ Y_{n} \} _{n\in
\mathbb{T}}$ which, conditionally on~$ \{ X_{n} \} $, are
independent and distributed according to
%
%e4.1 ###
%
\begin{equation}\label{eq:appobs}
Y_{n} \vert\{X_{n}=x\} \sim g ( x,\cdot) .
\end{equation}
This family of models is important, because almost every stationary time
series model appearing in the literature can be cast into this form.
Given $%
y_{1\dvtx P}$, we are often interested in computing the sequence of posterior\vadjust{\goodbreak}
distributions $ \{ p ( \mathbf{x}_{n} \vert
y_{1\dvtx n} ) \} _{n\in\mathbb{T}}$ to perform goodness-of-fit
and/or to compute the marginal likelihood $p ( y_{1\dvtx P} ) $. By
defining the unnormalized distribution as
%
%e4.2 ###
%
\begin{equation}\label{eq:filtering}
\gamma_{n} ( \mathbf{x}_{n} ) =p ( \mathbf{x}%
_{n},y_{1\dvtx n} ) =\mu( x_{1} ) g ( x_{1},y_{1} )
\prod_{k=2}^{n}f ( x_{k-1},x_{k} ) g ( x_{k},y_{k} )
\end{equation}
(which is typically known pointwise), we have $\pi_{n} ( \mathbf{x}%
_{n} ) =p ( \mathbf{x}_{n} \vert y_{1\dvtx n} ) $ and $%
Z_{n}=p ( y_{1\dvtx n} ) $ so that SIMCMC can be applied.

We will consider here three examples where the SIMCMC algorithms are
compared to their SMC counterparts. For a fixed number of
iterations/particles, SMC and SIMCMC have approximately the same
computational complexity. The same proposals and the same number of samples
were thus used to allow for a fair comparison. Note that we chose not
to use
any burn-in period for the SIMCMC and we initialize the algorithm by picking
$\mathbf{x}_{n}^{ ( 0 ) }= ( \mathbf{x}_{n-1}^{ ( 0 )
},x_{n}^{ ( 0 ) } ) $ for any $n$ where $\mathbf{x}%
_{P}^{ ( 0 ) }$ is a sample from the prior. The SMC algorithms
were implemented using a stratified resampling procedure~\cite{kitagawa1996}. The first two examples compare SMC to SIMCMC in
terms of log-likelihood
estimates. The third example demonstrates the use of SIMCMC in a real-time
tracking application.

%s4.1 ###
\subsection{Linear Gaussian model}

We consider a linear Gaussian model where $E=\mathbb{R}^{d}$:
%
%e4.3 ###
%
\begin{eqnarray} \label{eq:evollineargauss}
X_{n}& = & AX_{n-1}+\sigma_{v}V_{n}, \nonumber\\[-8pt]\\[-8pt]
Y_{n}& = & X_{n}+\sigma_{w}W_{n}, \nonumber
\end{eqnarray}
with $X_{1}\sim\mathcal{N} ( 0,\Delta)$, $V_{n}
\stackrel{\mathrm{i.i.d.}}{\sim}\mathcal{N} ( 0,\Delta) $,
$W_{n}\stackrel{\mathrm{i.i.d.}}{\sim}\mathcal{N} ( 0,\Delta) $,
$\Delta=\operatorname{diag}(
1,\ldots,1 ) $ and $A$ is a random (doubly stochastic) matrix. Here $%
\mathcal{N} ( \mu,\Sigma) $ is a Gaussian distribution of mean $%
\mu$ and variance-covariance matrix $\Sigma$. For this model we can
compute the marginal likelihood $Z_{P}=p ( y_{1\dvtx P} ) $ exactly
using the Kalman filter. This allows us to compare our results to the ground
truth.

We use two proposal densities: the prior density $f (
x_{n-1},x_{n} ) $ and the optimal density (\ref{eq:evollineargauss})
given by $q_{n}^{\mathrm{opt}} ( \mathbf{x}_{n-1},x_{n} ) \propto
f ( x_{n-1},x_{n} ) g ( x_{n},y_{n} ) $ which is a
Gaussian. In both cases, it is easy to check that Assumption \ref
{assumA1} is satisfied.

For $d=2,5,10$, we simulated a realization of $P=100$ observations for $
\sigma_{v}=2$ and $\sigma_{w}=0.5$. In Tables~\ref{table1} and \ref
{table2}, we
present the
performance of both SIMCMC and SMC for a varying $d$, a varying number of
samples and the two proposal distributions in terms on Root Mean Square
Error (RMSE) of the log-likelihood estimate.

%t1
%
\begin{table}
\caption{RMSE for SMC and SIMCMC algorithms over 100 realizations for
prior proposal}
\label{table1}
\begin{tabular*}{\tablewidth}{@{\extracolsep{\fill
}}ld{2.2}d{2.2}d{2.2}d{2.2}d{2.2}@{}}
\hline
$\bolds N$ & \multicolumn{1}{c}{\textbf{1000}} &
\multicolumn{1}{c}{\textbf{2500}}
& \multicolumn{1}{c}{\textbf{5000}} & \multicolumn{1}{c}{\textbf{10,000}}
& \multicolumn{1}{c@{}}{\textbf{25,000}} \\
\hline
SMC, $d=2$ & 1.66 & 0.98 & 0.63 & 0.52 & 0.29 \\
SIMCMC, $d=2$ & 1.57 & 0.97 & 0.75 & 0.59 & 0.41 \\
SMC, $d=5$ & 4.84 & 4.76 & 3.06 & 2.18 & 1.59 \\
SIMCMC, $d=5$ & 5.57 & 5.41 & 4.12 & 2.36 & 1.83 \\
SMC, $d=10$ & 16.91 & 14.57 & 11.14 & 10.61 & 8.91 \\
SIMCMC, $d=10$ & 18.22 & 16.78 & 14.56 & 12.46 & 11.25 \\
\hline
\end{tabular*}
\end{table}

%t2
%
\begin{table}[b]
\caption{RMSE for SMC and SIMCMC algorithms over 100 realizations for
optimal proposal}
\label{table2}
\begin{tabular*}{\tablewidth}{@{\extracolsep{\fill}}lccccc@{}}
\hline
$\bolds N$ & \multicolumn{1}{c}{\textbf{1000}} &
\multicolumn{1}{c}{\textbf{2500}}
& \multicolumn{1}{c}{\textbf{5000}} & \multicolumn{1}{c}{\textbf{10,000}}
& \multicolumn{1}{c@{}}{\textbf{25,000}} \\
\hline
SMC, $d=2$ & 0.33 & 0.17 & 0.09 & 0.06 & 0.04 \\
SIMCMC, $d=2$ & 0.37 & 0.19 & 0.14 & 0.11 & 0.06 \\
SMC, $d=5$ & 0.28 & 0.16 & 0.10 & 0.07 & 0.06 \\
SIMCMC, $d=5$ & 0.29 & 0.23 & 0.15 & 0.12 & 0.07 \\
SMC, $d=10$ & 0.18 & 0.14 & 0.09 & 0.05 & 0.07 \\
SIMCMC, $d=10$ & 0.31 & 0.20 & 0.16 & 0.12 & 0.10 \\
\hline
\end{tabular*}
\end{table}

As expected, the performance of our estimates is very significantly improved
when the optimal distribution is used as the observations are quite
informative. Unsurprisingly, SMC outperform SIMCMC for a fixed
computational complexity.

%s4.2 ###
\subsection{A nonlinear non-Gaussian state-space model}

We now consider a nonlinear non-Gaussian state-space model introduced in
\cite{kitagawa1996} which has been used in many SMC publications:
\begin{eqnarray*}
X_{n}& = &\frac{X_{n-1}}{2}+\frac{25X_{n-1}}{1+X_{n-1}^{2}}+8\cos(
1.2n ) +V_{n}, \\
Y_{n}& = &\frac{X_{n}^{2}}{20}+W_{n},
\end{eqnarray*}
where $X_{1}\sim\mathcal{N} ( 0,5 )$, $V_{n}\stackrel{\mathrm{i.i.d.}%
}{\sim}\mathcal{N} ( 0,\sigma_{v}^{2} ) $ and $W_{n}\stackrel{%
\mathrm{i.i.d.}}{\sim}\mathcal{N} ( 0,\sigma_{w}^{2} ) $. As the
sign of the state $X_{n}$ is not observed, the posterior distributions $
\{ p ( x_{1\dvtx n} \vert y_{1\dvtx n} ) \} $ are
multimodal. SMC approximations are able to capture properly the
multimodality of the posteriors. This allows us to assess here whether
SIMCMC can also explore properly these multimodal distributions by comparing
SIMCMC estimates to SMC estimates.

We use as a proposal density the prior density $f ( x_{n-1},x_{n} )
$. In this case, it is easy to check that Assumption~\ref{assumA1} is satisfied.

%t3
%
\begin{table}
\caption{Average RMSE of log-likelihood estimates for SMC and SIMCMC
algorithms over 100 realizations}
\label{table3}
\begin{tabular*}{\tablewidth}{@{\extracolsep{\fill}}lccccc@{}}
\hline
$\bolds N$ & \multicolumn{1}{c}{\textbf{2500}} & \multicolumn
{1}{c}{\textbf{5000}}
& \multicolumn{1}{c}{\textbf{10,000}} & \multicolumn{1}{c}{\textbf{25,000}}
& \multicolumn{1}{c@{}}{\textbf{50,000}} \\
\hline
SMC, $\sigma_{w}^{2}=1$ & 0.80 & 0.55 & 0.40 & 0.24 & 0.17 \\
SIMCMC, $\sigma_{w}^{2}=1$ & 0.95 & 0.60 & 0.75 & 0.59 & 0.41 \\
SMC, $\sigma_{w}^{2}=2$ & 0.33 & 0.23 & 0.17 & 0.11 & 0.07 \\
SIMCMC, $\sigma_{w}^{2}=2$ & 0.91 & 0.70 & 0.50 & 0.38 & 0.29 \\
SMC, $\sigma_{w}^{2}=5$ & 0.13 & 0.10 & 0.08 & 0.05 & 0.03 \\
SIMCMC, $\sigma_{w}^{2}=5$ & 0.28 & 0.21 & 0.19 & 0.12 & 0.08 \\
\hline
\end{tabular*}
\end{table}

In Table~\ref{table3}, we present the performance of both SIMCMC and
SMC for a varying
number of samples and a varying $\sigma_{w}^{2}$ whereas we set $\sigma
_{v}^{2}=5$. Both SMC and SIMCMC are performing better as the signal to
noise ratio degrades. This should not come as a surprise. As we are using
the prior as a proposal, it is preferable to have a diffuse likelihood for
good performance. Experimentally we observed that SIMCMC and SMC estimates
coincide for large $N$. However for a fixed computational complexity, SIMCMC
is outperformed by SMC in terms of RMSE.

%s4.3 ###
\subsection{Target tracking}

We consider here a target tracking problem~\cite{Gilks99,Liu2001}.
The target is modeled using a standard constant velocity model
%
%e4.4 ###
%
\begin{equation} \label{eq:evoltarget}
X_{n}= \pmatrix{
1 & T & 0 & 0 \cr
0 & 1 & 0 & 0 \cr
0 & 0 & 1 & T \cr
0 & 0 & 0 & 1}
X_{n-1}+V_{n},
\end{equation}
where $V_{n}\stackrel{\mathrm{i.i.d.}}{\sim}\mathcal{N} ( 0,\Sigma)
$, with $T=1$ and
\[
\Sigma=5 \pmatrix{
T^{3}/3 & T^{2}/2 & 0 & 0 \cr
T^{2}/2 & T & 0 & 0 \cr
0 & 0 & T^{3}/3 & T^{2}/2 \cr
0 & 0 & T^{2}/2 & T}.
\]
The state vector $X_{n}= (
X_{n}^{1},X_{n}^{2},X_{n}^{3},X_{n}^{4} ) ^{\mathrm{T}}$ is such that $%
X_{n}^{1}$ (resp., $X_{n}^{3}$) corresponds to the horizontal (resp.,
vertical) position of the target whereas $X_{n}^{2}$ (resp., $X_{n}^{4}$)
corresponds to the horizontal (resp., vertical) velocity. In many real
tracking applications, the observations are collected at random times~\cite{fox2001}. We have the following measurement equation:
%
%e4.5 ###
%
\begin{equation}\label{eq:obstarget}
Y_{n}=\tan^{-1} \biggl( \frac{X_{n}^{3}}{X_{n}^{1}} \biggr) +W_{n},
\end{equation}
where $W_{n}\stackrel{\mathrm{i.i.d.}}{\sim}\mathcal{N} ( 0,10^{-2} )
$; these parameters are representative of real-world tracking
scenarios. We
assume that we only observe $Y_{n}$ at the time indexes $n=4k$ where
$k\in
\mathbb{N}$ and, when $n\neq4k$, we observe $Y_{n}$ with probability $%
p=0.25 $. We are here interested in estimating the associated sequence of
posterior distributions on-line.

Assume the computational complexity of the SMC method is such that only
$%
N=1000$ particles can be used in one time unit. In such scenarios, we can
either use SMC with $N$ particles to estimate the sequence of posterior
distributions of interest or SMC with say $N^{\prime}=4000$ particles and
chose to ignore observations that might appear when $n\neq4k$. These are
two standard approaches used in applications. Alternatively, we can use the
SIMCMC method to select adaptively the number of particles as discussed in
Section~\ref{subsec:extensions}. If SIMCMC algorithm only adds particles
to the approximation of the current posterior at time $n$, it will use
approximately $mN$ particles to approximate this posterior if the next
observation only appears at time $n+m$.

%t4
%
\begin{table}
\caption{Average RMSE of the Monte Carlo state estimate}
\label{table4}
\begin{tabular*}{\tablewidth}{@{\extracolsep{\fill}}lccc@{}}
\hline
\textbf{Algorithm} & \textbf{SMC with} $\bolds N$
& \textbf{SMC with} $\bolds{N^{\prime}}$ & \textbf{SIMCMC} \\
\hline
Average RMSE & 2.14 & 3.21 & 1.62 \\
\hline
\end{tabular*}
\end{table}

We simulated 50 realizations of $P=100$ observations according to the
model (\ref{eq:evoltarget}) and (\ref{eq:obstarget}) and use as a
proposal density the
prior density $f ( x_{n-1},\break x_{n} ) $. This ensures that Assumption
\ref{assumA1} is satisfied. In Table~\ref{table4}, we display the
performance of SMC with $N$
particles, $N^{\prime}$ particles (ignoring some observations) and SIMCMC
using an adaptive number of particles in terms of the average RMSE of the
estimate of the conditional expectation $\mathbb{E} (
X_{n} \vert y_{1\dvtx n} ) $. In such scenarios, SIMCMC methods clearly
outperforms SMC methods.

%s5 ###
\section{Discussion}\label{sec:discussion}

We have described an iterative algorithm based on interacting non-Markovian
sequences which is an alternative to SMC and have established convergence
results validating this methodology. The algorithm is straightforward to
implement for people already familiar with MCMC. The main strength of
\mbox{SIMCMC}
compared to SMC is that it allows us to iteratively improve our estimates
in an MCMC-like fashion until a suitable stopping criterion is met.
This is
useful as in numerous applications the number of particles required to
ensure the estimates are reasonably precise is unknown. It is also
useful in
real-time applications when one is unsure of exactly how much time will be
available between successive arrivals of data points.

As discussed in Section~\ref{subsec:extensions}, numerous variants of
SIMCMC can be easily developed which have no SMC equivalent. The fact that
such schemes do not need to define a target distribution on an extended
state-space admitting $\{\pi_{n}\}_{n\in\mathbb{T}}$ as marginals offers
a lot of flexibility. For example, if we have access to multiple processors,
it is possible to sample from each $\pi_{n}$ independently using standard
MCMC and perform several interactions simultaneously. Adaptive versions of
the algorithms can also be proposed by monitoring the acceptance ratio of
the MH steps. If the acceptance probability of the MH move between say
$\pi
_{n-1}$ and $\pi_{n}$ is low, we could, for example, increase the
number of
proposals at this time index.

From a theoretical point of view, there are a number of interesting
questions to explore. Under additional stability assumptions on the
Feynman--Kac semigroup induced by $\{\pi_{n}\}_{n\in\mathbb{T}}$ and $%
\{q_{n}\}_{n\in\mathbb{T}}$~\cite{Delm04}, Chapter 4, we have recently
established in~\cite{DelMoral2009} convergence results ensuring that, for
functions of the form $f_{n} ( \mathbf{x}_{n} ) =f_{n} (
x_{n} ) $, the bound $C_{1,n}$ in Theorem \ref%
{theorem:convergenceofaverages} is independent of $n$. A central limit
theorem has also been established in~\cite{bercu2008}.

\begin{appendix}\label{app}
%s6 ###
\section*{Appendix: Proofs of convergence}

Our proofs rely on the Poisson equation~\cite{glynn1996} and are
inspired by
ideas developed in~\cite{andrieumoulines2006,andrieuajay2007,delmoralmiclo2003}.
However, contrary to standard adaptive MCMC schemes
\cite{andrieumoulines2006,robertsrosenthal2007}, the Markov kernels
we consider do not necessarily admit the target distribution as invariant
distribution; see~\cite{delmoralmiclo2003} for similar scenarios. However,
in our context, it is still possible to establish stronger results than in
\cite{delmoralmiclo2003} as we can characterize exactly these invariant
distributions; see Proposition~\ref{lemma:invariantdistributionperturbedkernel}.

%s6.1 ###
\subsection{Notation}

We denote by $\mathcal{P} ( E_{n} ) $ the set of probability
measures on $ ( E_{n},\mathcal{F}_{n} ) $. We introduce the
independent Metropolis--Hastings (MH) kernel $K_{1}\dvtx
E_{1}\times\mathcal
{F}%
_{1}\rightarrow[ 0,1 ] $ defined by
%
%e6.1 ###
%
\begin{eqnarray} \label{eq:initialMHkernel}
K_{1} ( \mathbf{x}_{1},d\mathbf{x}_{1}^{\prime} ) & = &\alpha
_{1} ( \mathbf{x}_{1},\mathbf{x}_{1}^{\prime} ) q_{1} ( d%
\mathbf{x}_{1}^{\prime} ) \nonumber\\[-8pt]\\[-8pt]
&&{} + \biggl( 1-\int\alpha_{1} ( \mathbf{x}_{1},\mathbf{y}_{1} )
q_{1} ( d\mathbf{y}_{1} ) \biggr) \delta_{\mathbf{x}_{1}} ( d%
\mathbf{x}_{1}^{\prime} ) . \nonumber
\end{eqnarray}
For $n\in\{ 2,\ldots,P \} $, we associate with any $\mu_{n-1}\in
\mathcal{P} ( E_{n-1} ) $ the Markov kernel $K_{n,\mu
_{n-1}}\dvtx E_{n}\times\mathcal{F}_{n}\rightarrow[ 0,1 ] $ defined
by
%
%e6.2 ###
%
\begin{eqnarray}\label{eq:perturbedMHkernel}\qquad
K_{n,\mu_{n-1}} ( \mathbf{x}_{n},d\mathbf{x}_{n}^{\prime} )
&=&\alpha_{n} ( \mathbf{x}_{n},\mathbf{x}_{n}^{\prime} ) ( \mu
_{n-1}\times q_{n} ) ( d\mathbf{x}_{n}^{\prime} )
\nonumber\\[-8pt]\\[-8pt]
&&{} + \biggl( 1-\int\alpha_{n} ( \mathbf{x}_{n},\mathbf{y}_{n} )
( \mu_{n-1}\times q_{n} ) ( d\mathbf{y}_{n} ) \biggr)
\delta_{\mathbf{x}_{n}} ( d\mathbf{x}_{n}^{\prime} ), \nonumber
\end{eqnarray}
where $\mathbf{x}_{n}^{\prime}= ( \mathbf{x}_{n-1}^{\prime
},x_{n}^{\prime} ) $. In (\ref{eq:initialMHkernel}) and (\ref%
{eq:perturbedMHkernel}), we have for $n\in\mathbb{T}$
\[
\alpha_{n} ( \mathbf{x}_{n},\mathbf{x}_{n}^{\prime} ) =1\wedge
\frac{w_{n} ( \mathbf{x}_{n}^{\prime} ) }{w_{n} ( \mathbf{x}%
_{n} ) }.
\]

We use $ \Vert\cdot\Vert_{\mathrm{tv}}$ to denote the
total variation norm and for any Markov kernel
\[
K^{i} ( \mathbf{x},d\mathbf{x}^{\prime} ) :=\int K^{i-1} (
\mathbf{x},d\mathbf{y} ) K ( \mathbf{y},d\mathbf{x}^{\prime
} ) .
\]

%s6.2 ###
\subsection{Preliminary results}

We establish here the expression of the invariant distributions of the
kernels $K_{1} ( \mathbf{x}_{1},d\mathbf{x}_{1}^{\prime} ) $ and $%
K_{n,\mu_{n-1}} ( \mathbf{x}_{n},d\mathbf{x}_{n}^{\prime} ) $ and
establish that these kernels are geometrically ergodic. We also provide some
perturbation bounds for $K_{n,\mu_{n-1}} ( \mathbf{x}_{n},d\mathbf{x}%
_{n}^{\prime} ) $ and its invariant distribution.
\begin{lemma}
\label{lemma:k1geometricallyergodic}
Given Assumption~\ref{assumA1}, $K_{1} ( \mathbf{x}_{1},d
\mathbf{x}_{1}^{\prime} ) $ is uniformly geometrically ergodic of
invariant distribution $\pi_{1} ( d\mathbf{x}_{1} ) $.
\end{lemma}

By construction, $K_{1} ( \mathbf{x}_{1},d\mathbf{x}_{1}^{\prime
} ) $ is an MH algorithm of invariant distribution $\pi_{1} ( d%
\mathbf{x}_{1} ) $. Uniform ergodicity follows from Assumption \ref
{assumA1}; see, for example,
Theorem 2.1. in~\cite{mengersentweedie1996}.
\begin{proposition}
\label{lemma:invariantdistributionperturbedkernel}Given Assumption \ref
{assumA1}, for any $%
n\in\{ 2,\ldots,P \} $ and any $\mu_{n-1}\in\mathcal{P} (
E_{n-1} ) $, $K_{n,\mu_{n-1}} ( \mathbf{x}_{n},d\mathbf{x}%
_{n}^{\prime} ) $ is uniformly geometrically ergodic of invariant
distribution
%
%e6.3 ###
%
\begin{equation} \label{eq:invariantdistributionperturbedkernel}
\omega_{n} ( \mu_{n-1} ) ( d\mathbf{x}_{n} ) =\frac{%
\pi_{n/n-1} ( \mathbf{x}_{n-1} ) \cdot ( \mu_{n-1}\times\overline{%
\pi}_{n} ) ( d\mathbf{x}_{n} ) }{\mu_{n-1} (
\pi_{n/n-1} ) },
\end{equation}
where $\overline{\pi}_{n} ( \mathbf{x}_{n-1},dx_{n} ) $ and $%
\pi_{n/n-1} ( \mathbf{x}_{n-1} ) $ are defined, respectively, in (\ref%
{eq:conditionaldensitypin}) and~(\ref{eq:optimalimportanceweight}).
\end{proposition}
\begin{pf}
To establish the result, it is sufficient to rewrite
\begin{eqnarray*}
w_{n} ( \mathbf{x}_{n} ) & = & \frac{Z_{n}}{Z_{n-1}}\frac{{\pi
_{n} ( \mathbf{x}_{n} ) }/{\pi_{n-1} ( \mathbf{x}_{n-1} ) }%
\mu_{n-1} ( \mathbf{x}_{n-1} ) }{ ( \mu_{n-1}\times
q_{n} ) ( \mathbf{x}_{n} ) } \\
& = & \frac{Z_{n}}{Z_{n-1}}\frac{\pi_{n/n-1} ( \mathbf{x}_{n-1} )
( \mu_{n-1}\times\overline{\pi}_{n} ) ( \mathbf{x}%
_{n} ) }{ ( \mu_{n-1}\times q_{n} ) ( \mathbf{x}%
_{n} ) }.
\end{eqnarray*}
This shows that $K_{n,\mu_{n-1}} ( \mathbf{x}_{n},d\mathbf{x}%
_{n}^{\prime} ) $ is nothing but a standard MH algorithm of proposal
distribution $ ( \mu_{n-1}\times q_{n} ) ( \mathbf{x}%
_{n} ) $ and target distribution given by~(\ref%
{eq:invariantdistributionperturbedkernel}). This distribution is always
proper because Assumption~\ref{assumA1} implies that $\pi_{n/n-1} (
\mathbf{x}%
_{n-1} ) <\infty$ over $E_{n-1}$. Uniform ergodicity follows from
Theorem~2.1. in~\cite{mengersentweedie1996}.
\end{pf}
\begin{Corollary*} It follows that for any $n\in\{ 2,\ldots,P \} $
there exists $\rho_{n}<1$ such that for any $m\in\mathbb{N}_{0}$ and $%
\mathbf{x}_{n}\in E_{n}$
%
%e6.4 ###
%
\begin{equation} \label{eq:uniformergodicity}
\Vert K_{n,\mu_{n-1}}^{m} ( \mathbf{x}_{n},\cdot) -\omega
_{n} ( \mu_{n-1} ) ( \cdot) \Vert_{\mathrm{tv}}
\leq\rho_{n}^{m}.
\end{equation}
\end{Corollary*}
\begin{proposition}
\label{proposition:sensitivity}Given Assumption~\ref{assumA1}, for any
$n\in\{
2,\ldots,P \} $, we have for any $\mu_{n-1}$, $\nu_{n-1}\in\mathcal{P}%
( E_{n-1} ) $, $\mathbf{x}_{n}\in E_{n}$ and $m\in\mathbb{N}$
%
%e6.5 ###
%
\begin{equation} \label{eq:sensitivitykernel}
\Vert K_{n,\mu_{n-1}}^{m} ( \mathbf{x}_{n},\cdot)
-K_{n,\nu_{n-1}}^{m} ( \mathbf{x}_{n},\cdot) \Vert
_{\mathrm{tv}}\leq\frac{2}{1-\rho_{n}} \Vert\mu_{n-1}-\nu
_{n-1} \Vert_{\mathrm{tv}}
\end{equation}
and
%
%e6.6 ###
%
\begin{equation}\label{eq:sensitivityinvariantdistribution}
\Vert\omega_{n} ( \mu_{n-1} ) -\omega_{n} ( \nu
_{n-1} ) \Vert_{\mathrm{tv}}\leq\frac{2}{1-\rho_{n}}%
\Vert\mu_{n-1}-\nu_{n-1} \Vert_{\mathrm{tv}}.\vadjust{\goodbreak}
\end{equation}
\end{proposition}
\begin{pf}
For any $f_{n}\in B ( E_{n} ) $, we have the
following decomposition:
\begin{eqnarray*}
&& K_{n,\mu_{n-1}}^{m} ( f_{n} ) ( \mathbf{x}_{n} )
-K_{n,\nu_{n-1}}^{m} ( f_{n} ) ( \mathbf{x}_{n} ) \\
&&\qquad =\sum_{j=0}^{m-1}K_{n,\mu_{n-1}}^{j} ( K_{n,\mu_{n-1}}-K_{n,\nu
_{n-1}} ) K_{n,\nu_{n-1}}^{m-j-1} ( f_{n} ) ( \mathbf{x}%
_{n} ) \\
&&\qquad =\sum_{j=0}^{m-1}K_{n,\mu_{n-1}}^{j} ( K_{n,\mu_{n-1}}-K_{n,\nu
_{n-1}} ) \bigl( K_{n,\nu_{n-1}}^{m-j-1} ( f_{n} ) (
\mathbf{x}_{n} ) -\omega_{n} ( \nu_{n-1} ) (
f_{n} ) \bigr) .
\end{eqnarray*}
From Assumption~\ref{assumA1}, we know that for any $\nu_{n-1}\in
\mathcal{P} (
E_{n-1} ) $
\[
\Vert K_{n,\nu_{n-1}}^{m-j-1} ( \mathbf{x}_{n},\cdot)
-\omega_{n} ( \nu_{n-1} ) \Vert_{\mathrm{tv}}\leq
\rho_{n}^{m-j-1}
\]
and from (\ref{eq:perturbedMHkernel}) for any $\mathbf{x}_{n}\in E_{n}$
and $%
f_{n}\in B ( E_{n} ) $
\begin{eqnarray*}
&& K_{n,\mu_{n-1}} ( f_{n} ) ( \mathbf{x}_{n} ) -K_{n,\nu
_{n-1}} ( f_{n} ) ( \mathbf{x}_{n} ) \\
&&\qquad =\int f_{n} ( \mathbf{x}_{n}^{\prime} ) \alpha_{n} (
\mathbf{x}_{n},\mathbf{x}_{n}^{\prime} ) \bigl( ( \mu_{n-1}-\nu
_{n-1} ) \times q_{n} \bigr) ( d\mathbf{x}_{n}^{\prime} ) \\
&&\qquad\quad{} +f_{n} ( \mathbf{x}_{n} ) \int\alpha_{n} ( \mathbf{x}_{n},%
\mathbf{y}_{n}^{\prime} ) \bigl( ( \nu_{n-1}-\mu_{n-1} )
\times q_{n} \bigr) ( d\mathbf{y}_{n}^{\prime} )
\end{eqnarray*}
thus%
\[
\Vert K_{n,\mu_{n-1}} ( \mathbf{x}_{n},\cdot) -K_{n,\nu
_{n-1}} ( \mathbf{x}_{n},\cdot) \Vert_{\mathrm{tv}}%
\leq2 \Vert\mu_{n-1}-\nu_{n-1} \Vert_{\mathrm{tv}}.
\]
So
\begin{eqnarray*}
\Vert K_{n,\mu_{n-1}}^{m} ( \mathbf{x}_{n},\cdot) -K_{n,\nu
_{n-1}}^{m} ( \mathbf{x}_{n},\cdot) \Vert_{\mathrm{tv}}
& \leq & 2 \Vert\mu_{n-1}-\nu_{n-1} \Vert_{\mathrm{tv}}%
\sum_{j=0}^{m-1}\rho_{n}^{m-j-1} \\
& = & 2\frac{1-\rho_{n}^{m}}{1-\rho_{n}} \Vert\mu_{n-1}-\nu
_{n-1} \Vert_{\mathrm{tv}}.
\end{eqnarray*}
Hence, (\ref{eq:sensitivitykernel}) follows and we obtain (\ref%
{eq:sensitivityinvariantdistribution}) by taking the limit as
$m\rightarrow
\infty$.
\end{pf}

%s6.3 ###
\subsection{Convergence of averages}

For any $n\in\{ 2,\ldots,P \} $, $p\geq1$ and $f_{n}\in B (
E_{n} ) $ we want to study%
\[
\mathbb{E}_{\mathbf{x}_{1\dvtx n}^{ ( 0 ) }} \bigl[ \bigl\vert\widehat{%
\pi}{}^{ ( i ) }_{n} ( f_{n} ) -\pi_{n} ( f_{n} )
\bigr\vert^{p} \bigr] ^{1/p}.
\]
We have
%
%e6.7 ###
%
\begin{equation} \label{eq:decompositionaverages}
\widehat{\pi}{}^{ ( i ) }_{n} ( f_{n} ) -\pi_{n} (
f_{n} ) =\widehat{\pi}{}^{ ( i ) }_{n} ( f_{n} )
-S_{n}^{ ( i ) } ( f_{n} ) +S_{n}^{ ( i ) } (
f_{n} ) -\pi_{n} ( f_{n} ),
\end{equation}
where
\[
S_{n}^{ ( i ) } ( f_{n} ) =\frac{1}{i+1}%
\sum_{m=0}^{i}\omega_{n} \bigl( \widehat{\pi}{}^{ ( m )
} _{n-1}\bigr) ( f_{n} ) .
\]
To study the first term on the right-hand side of (\ref{eq:decompositionaverages}), we
introduce the Poisson equation~\cite{glynn1996}
\[
f_{n} ( x ) -\omega_{n} ( \mu) ( f_{n} ) =%
\widehat{f}_{n,\mu} ( x ) -K_{n,\mu} ( \widehat{f}_{n,\mu
} ) ( x ),
\]
whose solution, if $K_{n,\mu}$ is geometrically ergodic, is given by
%
%e6.8 ###
%
\begin{equation} \label{eq:solutionpoissonequation}
\widehat{f}_{n,\mu} ( x ) =\sum_{i\in\mathbb{N}_{0}} [
K_{n,\mu}^{i} ( f_{n} ) ( x ) -\omega_{n} ( \mu
) ( f_{n} ) ] .
\end{equation}
We have
%
%e6.9 ###
%
\begin{eqnarray} \label{eq:longdecompositionaverages}
&&( i+1 ) \bigl[ \widehat{\pi}_{n}^{ ( i ) } (
f_{n} ) -S_{n}^{ ( i ) } \bigr] \nonumber\\
&&\qquad =  M_{n}^{ ( i+1 )
} ( f_{n} )
+\sum_{m=0}^{i} \bigl[ \widehat{f}_{n,\widehat{\pi}_{n-1}^{ (
m+1 ) }} \bigl( \mathbf{X}_{n}^{ ( m+1 ) } \bigr) -\widehat{f}%
_{n,\widehat{\pi}_{n-1}^{ ( m ) }} \bigl( \mathbf{X}_{n}^{ (
m+1 ) } \bigr) \bigr] \\
&&\qquad\quad{} +\widehat{f}_{n,\widehat{\pi}_{n-1}^{ ( 0 ) }} \bigl( \mathbf{X}%
_{n}^{ ( 0 ) } \bigr) -\widehat{f}_{n,\widehat{\pi}_{n-1}^{ (
i+1 ) }} \bigl( \mathbf{X}_{n}^{ ( i+1 ) } \bigr), \nonumber
\end{eqnarray}
where
%
%e6.10 ###
%
\begin{equation} \label{eq:martingaleterm}
M_{n}^{ ( i ) } ( f_{n} ) =\sum_{m=0}^{i-1} \bigl[ \widehat{%
f}_{n,\widehat{\pi}_{n-1}^{ ( m ) }} \bigl( \mathbf{X}_{n}^{ (
m+1 ) } \bigr) -K_{\widehat{\pi}_{n-1}^{ ( m ) }} (
\widehat{f}_{n,\widehat{\pi}_{n-1}^{ ( m ) }} ) \bigl(
\mathbf{X}_{n}^{ ( m ) } \bigr) \bigr]
\end{equation}
is a $\mathcal{G}_{n}^{i}$-martingale with $\mathcal{G}_{n}^{i}=\sigma
( \mathbf{X}_{1}^{ ( 1\dvtx i ) },\mathbf{X}_{2}^{ ( 1\dvtx i ) },\ldots
,\mathbf{X}_{n}^{ ( 1\dvtx i ) } ) $ where we define
$\mathbf{X}_{k}^{ ( 1\dvtx i ) }= ( \mathbf{X}%
_{k}^{ ( 1 ) },\ldots,\mathbf{X}_{k}^{ ( i ) } ) $.

We remind the reader that $B ( E_{n} ) =\{
f_{n}\dvtx E_{n}\rightarrow\mathbb{R}$ such that $ \Vert
f_{n} \Vert\leq1\} $ where $ \Vert f_{n} \Vert=%
\sup_{\mathbf{x}_{n}\in E_{n}} \vert f_{n} ( \mathbf{x}%
_{n} ) \vert$. We establish the following propositions.
\begin{proposition}
\label{proposition:poissonbounded}Given Assumption~\ref{assumA1}, for
any $n\in\{
2,\ldots,P \} $, $\mathbf{x}_{1\dvtx n}^{ ( 0 ) }$, \mbox{$p\geq1$}, $%
f_{n}\in B ( E_{n} ) $ and $m\in\mathbb{N}_{0}$, we have%
\[
\mathbb{E}_{\mathbf{x}_{1\dvtx n}^{ ( 0 ) }} \bigl[ \bigl\vert\widehat{f}%
_{n,\widehat{\pi}_{n-1}^{ ( m ) }} \bigl( \mathbf{X}_{n}^{ (
m+1 ) } \bigr) \bigr\vert^{p} \bigr] ^{1/p}\leq\frac{1}{1-\rho_{n}}.
\]
\end{proposition}
\begin{pf}
Assumption~\ref{assumA1} ensures that $\widehat{f}_{n,\widehat{\pi}
_{n-1}^{ ( m ) }}$ is given by an expression of the form (\ref%
{eq:solutionpoissonequation}). We have
\begin{eqnarray*}
&& \mathbb{E}_{\mathbf{x}_{1\dvtx n}^{ ( 0 ) }} \bigl[ \bigl\vert\widehat{%
f}_{n,\widehat{\pi}_{n-1}^{ ( m ) }} \bigl( \mathbf{X}_{n}^{ (
m+1 ) } \bigr) \bigr\vert^{p} \bigr] ^{1/p} \\
&&\qquad \leq\sum_{i\in\mathbb{N}_{0}}\mathbb{E}_{\mathbf{x}_{1\dvtx n}^{ ( 0 )
}} \bigl[ \bigl\vert K_{n,\widehat{\pi}_{n-1}^{ ( m ) }}^{i} (
f_{n} ) \bigl( \mathbf{X}_{n}^{ ( m+1 ) } \bigr)
-\omega_{n} \bigl( \widehat{\pi}_{n-1}^{ ( m ) } \bigr) (
f_{n} ) \bigr\vert^{p} \bigr] ^{1/p} \\
&&\qquad \leq\sum_{i\in\mathbb{N}_{0}}\mathbb{E}_{\mathbf{x}_{1\dvtx n}^{ ( 0 )
}} \bigl[ \mathbb{E}_{\mathbf{x}_{1\dvtx n}^{ ( 0 ) }} \bigl(
\bigl\vert K_{n,\widehat{\pi}_{n-1}^{ ( m ) }}^{i} (
f_{n} ) \bigl( \mathbf{X}_{n}^{ ( m+1 ) } \bigr)
-\omega_{n} \bigl( \widehat{\pi}_{n-1}^{ ( m ) } \bigr) (
f_{n} ) \bigr\vert^{p} \vert\mathcal{G}_{n-1}^{m} \bigr) \bigr]
^{1/p} \\
&&\qquad
\leq\sum_{i\in\mathbb{N}_{0}}\rho_{n}^{i}=\frac{1}{1-\rho_{n}},
\end{eqnarray*}
using Minkowski's inequality and the fact that $K_{n,\widehat{\pi}%
_{n-1}^{ ( m ) }}$ is an uniformly ergodic Markov kernel
conditional upon $\mathcal{G}_{n-1}^{m}$ using Assumption~\ref{assumA1}.
\end{pf}
\begin{proposition}
\label{proposition:martingalebounded}Given Assumption~\ref{assumA1},
for any $n\in\{
2,\ldots,P \} $ and any $p\geq1$, there exist $B_{1,n},B_{2,p}<\infty$
such that for any $\mathbf{x}_{1\dvtx n}^{ ( 0 ) }$, $f_{n}\in B (
E_{n} ) $ and $m\in\mathbb{N}$
\[
\mathbb{E}_{\mathbf{x}_{1\dvtx n}^{ ( 0 ) }} \bigl[ \bigl\vert
M_{n}^{ ( m ) } ( f_{n} ) \bigr\vert^{p} \bigr]
^{1/p}\leq B_{1,n}B_{2,p} m^{{1/2}}.
\]
\end{proposition}
\begin{pf}
For $p>1$ we use Burkholder's inequality (e.g.,~\cite{shiryaev1996},
page 499); that is, there exist constants $C_{1,n},C_{2,p}<\infty$
such
that%
%
%e6.11 ###
%
\begin{eqnarray} \label{eq:martingaleMn}\quad
&& \mathbb{E}_{\mathbf{x}%
_{1\dvtx n}^{ ( 0 ) }} \bigl[ \bigl\vert M_{n}^{ ( m ) } (
f_{n} ) \bigr\vert^{p} \bigr] ^{1/p}\nonumber\\
&&\qquad\leq C_{1,n}C_{2,p}\mathbb{E}_{\mathbf{x}_{1\dvtx n}^{ ( 0 ) }} \Biggl[
\Biggl( \sum_{i=0}^{m-1} \bigl[ \widehat{f}_{n,\widehat{\pi}_{n-1}^{ (
i ) }} \bigl( \mathbf{X}_{n}^{ ( i+1 ) } \bigr)\\
&&\qquad\quad\hspace*{96.3pt}{}-K_{n,\widehat{%
\pi}_{n-1}^{ ( i ) }} ( \widehat{f}_{n,\widehat{\pi}%
_{n-1}^{ ( i ) }} ) \bigl( \mathbf{X}_{n}^{ ( i )
} \bigr) \bigr] ^{2} \Biggr) ^{p/2} \Biggr] ^{1/p}. \nonumber
\end{eqnarray}
For $p\in( 1,2 ) $, we can bound the left-hand side of (\ref{eq:martingaleMn}%
)%
\begin{eqnarray*}
&& \mathbb{E}_{\mathbf{x}_{1\dvtx n}^{ ( 0 ) }} \Biggl[ \Biggl(
\sum_{i=0}^{m-1} \bigl[ \widehat{f}_{n,\widehat{\pi}_{n-1}^{ ( i )
}} \bigl( \mathbf{X}_{n}^{ ( i+1 ) } \bigr) -K_{n,\widehat{\pi}%
_{n-1}^{ ( i ) }} ( \widehat{f}_{n,\widehat{\pi}_{n-1}^{ (
i ) }} ) \bigl( \mathbf{X}_{n}^{ ( i ) } \bigr) \bigr]
^{2} \Biggr) ^{p/2} \Biggr] ^{1/p} \\
&&\qquad \leq\mathbb{E}_{\mathbf{x}_{1\dvtx n}^{ ( 0 ) }} \Biggl[ \Biggl(
2\sum_{i=0}^{m-1} \bigl[ \bigl\vert\widehat{f}_{n,\widehat{\pi}%
_{n-1}^{ ( i ) }} \bigl( \mathbf{X}_{n}^{ ( i+1 ) } \bigr)
\bigr\vert^{2}+ \bigl\vert K_{n,\widehat{\pi}_{n-1}^{ ( i )
}} ( \widehat{f}_{n,\widehat{\pi}_{n-1}^{ ( i ) }} )
\bigl( \mathbf{X}_{n}^{ ( i ) } \bigr) \bigr\vert^{2} \bigr]
\Biggr) ^{p/2} \Biggr] ^{1/p} \\
&&\qquad \leq\mathbb{E}_{\mathbf{x}_{1\dvtx n}^{ ( 0 ) }} \Biggl[ \Biggl(
2\sum_{i=0}^{m-1} \bigl[ \bigl\vert\widehat{f}_{n,\widehat{\pi}%
_{n-1}^{ ( i ) }} \bigl( \mathbf{X}_{n}^{ ( i+1 ) } \bigr)
\bigr\vert^{2}+ \bigl\vert K_{n,\widehat{\pi}_{n-1}^{ ( i )
}} ( \widehat{f}_{n,\widehat{\pi}_{n-1}^{ ( i ) }} )
\bigl( \mathbf{X}_{n}^{ ( i ) } \bigr) \bigr\vert^{2} \bigr]
\Biggr) \Biggr] ^{1/2}
\end{eqnarray*}
using $ ( a-b ) ^{2}\leq2 ( a^{2}+b^{2} ) $ and Jensen's
inequality. Now using Jensen's inequality again, we have
\[
\mathbb{E}_{\mathbf{x}_{1\dvtx n}^{ ( 0 ) }} \bigl[ \bigl\vert K_{n,%
\widehat{\pi}_{n-1}^{ ( i ) }} ( \widehat{f}_{n,\widehat{\pi}%
_{n-1}^{ ( i ) }} ) \bigl( \mathbf{X}_{n}^{ ( i )
} \bigr) \bigr\vert^{2} \bigr] \leq\mathbb{E}_{\mathbf{x}_{1\dvtx n}^{ (
0 ) }} \bigl[ K_{n,\widehat{\pi}_{n-1}^{ ( i ) }} (
\vert\widehat{f}_{n,\widehat{\pi}_{n-1}^{ ( i )
}} \vert^{2} ) \bigl( \mathbf{X}_{n}^{ ( i ) } \bigr) %
\bigr]
\]
and using Proposition~\ref{proposition:poissonbounded}, we obtain the bound
\[
\mathbb{E}_{\mathbf{x}_{1\dvtx n}^{ ( 0 ) }} \Biggl[ \Biggl(
\sum_{i=0}^{m-1} \bigl[ \widehat{f}_{n,\widehat{\pi}_{n-1}^{ ( i )
}} \bigl( \mathbf{X}_{n}^{ ( i+1 ) } \bigr) -K_{n,\widehat{\pi}%
_{n-1}^{ ( i ) }} ( \widehat{f}_{n,\widehat{\pi}_{n-1}^{ (
i ) }} ) \bigl( \mathbf{X}_{n}^{ ( i ) } \bigr) \bigr]
^{2} \Biggr) ^{p/2} \Biggr] ^{1/p}\leq\frac{2}{1-\rho_{n}}m^{1/2}.
\]

For $p\geq2$, we we can bound the left-hand side of (\ref{eq:martingaleMn}) through
Minkows\-ki's inequality
\begin{eqnarray*}
&&\mathbb{E}_{\mathbf{x}_{1\dvtx n}^{ ( 0 ) }} \bigl[ \bigl\vert
M_{n}^{ ( m ) } ( f_{n} ) \bigr\vert^{p} \bigr] ^{1/p} \\
&&\qquad\leq C_{1,n}C_{2,p} \Biggl( \sum_{i=0}^{m-1}\mathbb{E}_{\mathbf{x}%
_{1\dvtx n}^{ ( 0 ) }} \bigl[ \bigl\vert\widehat{f}_{n,\widehat{\pi}%
_{n-1}^{ ( i ) }} \bigl( \mathbf{X}_{n}^{ ( i+1 ) } \bigr)\\
&&\qquad\quad\hspace*{94.31pt}{}
-K_{n,\widehat{\pi}_{n-1}^{ ( i ) }} ( \widehat{f}_{n,\widehat{%
\pi}_{n-1}^{ ( i ) }} ) \bigl( \mathbf{X}_{n}^{ (
i ) } \bigr) \bigr\vert^{p} \bigr] ^{2/p} \Biggr) ^{1/2}.
\end{eqnarray*}
Using Minkowski's inequality again
\begin{eqnarray*}
&& \mathbb{E}_{\mathbf{x}_{1\dvtx n}^{ ( 0 ) }} \bigl[ \bigl\vert\widehat{%
f}_{n,\widehat{\pi}_{n-1}^{ ( i ) }} \bigl( \mathbf{X}_{n}^{ (
i+1 ) } \bigr) -K_{n,\widehat{\pi}_{n-1}^{ ( i ) }} (
\widehat{f}_{n,\widehat{\pi}_{n-1}^{ ( i ) }} ) \bigl(
\mathbf{X}_{n}^{ ( i ) } \bigr) \bigr\vert^{p} \bigr] \\
&&\qquad \leq\bigl( \mathbb{E}_{\mathbf{x}_{1\dvtx n}^{ ( 0 ) }} \bigl[
\bigl\vert\widehat{f}_{n,\widehat{\pi}_{n-1}^{ ( i ) }} \bigl(
\mathbf{X}_{n}^{ ( i+1 ) } \bigr) \bigr\vert^{p} \bigr] ^{1/p}+%
\mathbb{E}_{\mathbf{x}_{1\dvtx n}^{ ( 0 ) }} \bigl[ \bigl\vert K_{n,%
\widehat{\pi}_{n-1}^{ ( i ) }} ( \widehat{f}_{n,\widehat{\pi}%
_{n-1}^{ ( i ) }} ) \bigl( \mathbf{X}_{n}^{ ( i )
} \bigr) \bigr\vert^{p} \bigr] ^{1/p} \bigr) ^{p}.
\end{eqnarray*}
Now from Proposition~\ref{proposition:poissonbounded} and Jensen's
inequality, we can conclude for $p\geq1$. For $p=1$, we use Davis'
inequality (e.g.,~\cite{shiryaev1996}, page 499) to obtain the result using
similar arguments which are not repeated here.
\end{pf}
\begin{proposition}
\label{proposition:poissonincrements}Given Assumption~\ref{assumA1},
for any $n\in\{
2,\ldots,P \} $ and $p\geq1$ there exists $B_{n}<\infty$ such that for any
$\mathbf{x}_{1\dvtx n}^{ ( 0 ) }$, $f_{n}\in B ( E_{n} ) $ and $%
m\in\mathbb{N}_{0}$
\[
\mathbb{E}_{\mathbf{x}_{1\dvtx n}^{ ( 0 ) }} \bigl[ \bigl\vert\widehat{f}%
_{n,\widehat{\pi}_{n-1}^{ ( m+1 ) }} \bigl( \mathbf{X}_{n}^{ (
m+1 ) } \bigr) -\widehat{f}_{n,\widehat{\pi}_{n-1}^{ ( m )
}} \bigl( \mathbf{X}_{n}^{ ( m+1 ) } \bigr) \bigr\vert^{p} \bigr]
^{1/p}\leq\frac{B_{n}}{m+2}.
\]
\end{proposition}
\begin{pf}
Our proof is based on the following key decomposition
established in Lemma 3.4. in~\cite{delmoralmiclo2003}:
%
%e6.12 ###
%
\begin{eqnarray} \label{eq:identitypoissonequationdelmoralmiclo}\quad
&& \widehat{f}_{n,\widehat{\pi}_{n-1}^{ ( m+1 ) }} \bigl( \mathbf{X}%
_{n}^{ ( m+1 ) } \bigr) -\widehat{f}_{n,\widehat{\pi}%
_{n-1}^{ ( m ) }} \bigl( \mathbf{X}_{n}^{ ( m+1 ) } \bigr)
+\omega_{n} \bigl( \widehat{\pi}_{n-1}^{ ( m+1 ) } \bigr) (
\widehat{f}_{n,\widehat{\pi}_{n-1}^{ ( m ) }} ) \nonumber\\
&&\qquad =\sum_{i,j\in\mathbb{N}_{0}} \bigl( \delta_{\mathbf{X}_{n}^{ (
m+1 ) }}-\omega_{n} \bigl( \widehat{\pi}_{n-1}^{ ( m+1 )
} \bigr) \bigr) K_{n,\widehat{\pi}_{n-1}^{ ( m+1 ) }}^{i} (
K_{n,\widehat{\pi}_{n-1}^{ ( m+1 ) }}-K_{n,\widehat{\pi}%
_{n-1}^{ ( m ) }} ) \\
&&\qquad\quad\hspace*{22.46pt}{} \times K_{n,\widehat{\pi}_{n-1}^{ ( m ) }}^{j} \bigl(
f_{n}-\omega_{n} \bigl( \widehat{\pi}_{n-1}^{ ( m ) } \bigr)
( f_{n} ) \bigr).\nonumber
\end{eqnarray}
We have
%
%e6.13 ###
%
\begin{eqnarray} \label{eq:boundonfirstterm}
&& \bigl\vert\bigl( \delta_{\mathbf{X}_{n}^{ ( m+1 ) }}-\omega
_{n} \bigl( \widehat{\pi}_{n-1}^{ ( m+1 ) } \bigr) \bigr) K_{n,%
\widehat{\pi}_{n-1}^{ ( m+1 ) }}^{i} ( K_{n,\widehat{\pi}%
_{n-1}^{ ( m+1 ) }}-K_{n,\widehat{\pi}_{n-1}^{ ( m )
}} )\nonumber\hspace*{-25pt}\\
&&\quad\hspace*{72.4pt}{}\times K_{n,\widehat{\pi}_{n-1}^{ ( m ) }}^{j} \bigl(
f_{n}-\omega_{n} \bigl( \widehat{\pi}_{n-1}^{ ( m ) } \bigr)
( f_{n} ) \bigr) \bigr\vert\nonumber\hspace*{-25pt}\\
&&\qquad = \bigl\vert\bigl( \delta_{\mathbf{X}_{n}^{ ( m+1 ) }}-\omega
_{n} \bigl( \widehat{\pi}_{n-1}^{ ( m+1 ) } \bigr) \bigr) K_{n,%
\widehat{\pi}_{n-1}^{ ( m+1 ) }}^{i} ( K_{n,\widehat{\pi}%
_{n-1}^{ ( m+1 ) }}-K_{n,\widehat{\pi}_{n-1}^{ ( m )
}} ) K_{n,\widehat{\pi}_{n-1}^{ ( m ) }}^{j} (
f_{n} ) \bigr\vert\nonumber\hspace*{-25pt}\\
&&\qquad \leq\rho_{n}^{j} \bigl\Vert\bigl( \delta_{\mathbf{X}_{n}^{ (
m+1 ) }}-\omega_{n} \bigl( \widehat{\pi}_{n-1}^{ ( m+1 )
} \bigr) \bigr) K_{n,\widehat{\pi}_{n-1}^{ ( m+1 ) }}^{i} (
K_{n,\widehat{\pi}_{n-1}^{ ( m+1 ) }}-K_{n,\widehat{\pi}%
_{n-1}^{ ( m ) }} ) \bigr\Vert_{\mathrm{tv}} \hspace*{-25pt}\\
&&\qquad \leq\rho_{n}^{j}\times\frac{2}{1-\rho_{n}} \bigl\Vert\widehat{\pi}%
_{n-1}^{ ( m+1 ) }-\widehat{\pi}_{n-1}^{ ( m )
} \bigr\Vert_{\mathrm{tv}} \bigl\Vert\bigl( \delta_{\mathbf{X}%
_{n}^{ ( m+1 ) }}-\omega_{n} \bigl( \widehat{\pi}_{n-1}^{ (
m+1 ) } \bigr) \bigr) K_{n,\widehat{\pi}_{n-1}^{ ( m+1 )
}}^{i} \bigr\Vert_{\mathrm{tv}} \nonumber\hspace*{-25pt}\\
&&\qquad \leq\rho_{n}^{j}\times\frac{2}{1-\rho_{n}} \bigl\Vert\widehat{\pi}%
_{n-1}^{ ( m+1 ) }-\widehat{\pi}_{n-1}^{ ( m )
} \bigr\Vert_{\mathrm{tv}}\times\rho_{n}^{i} \bigl\Vert\delta_{%
\mathbf{X}_{n}^{ ( m+1 ) }}-\omega_{n} \bigl( \widehat{\pi}%
_{n-1}^{ ( m+1 ) } \bigr) \bigr\Vert_{\mathrm{tv}} \nonumber
\nonumber\hspace*{-25pt}\\
&&\qquad \leq\frac{2\rho_{n}^{i+j}}{1-\rho_{n}} \bigl\Vert\widehat{\pi}%
_{n-1}^{ ( m+1 ) }-\widehat{\pi}_{n-1}^{ ( m )
} \bigr\Vert_{\mathrm{tv}},\nonumber\hspace*{-25pt}
\end{eqnarray}
using Assumption~\ref{assumA1}, (\ref{eq:sensitivitykernel}) in Proposition
\ref{proposition:sensitivity} and Assumption~\ref{assumA1} again.

Now we have%
%
%e6.14 ###
%
\begin{eqnarray}\label{eq:boundon2ndterm}
&&\bigl\vert\omega_{n} \bigl( \widehat{\pi}_{n-1}^{ ( m+1 ) } \bigr)
( \widehat{f}_{n,\widehat{\pi}_{n-1}^{ ( m ) }} )\bigr\vert\nonumber\\
&&\qquad = \biggl\vert\omega_{n} \bigl( \widehat{\pi}_{n-1}^{ (
m+1 ) } \bigr) \biggl( \sum_{i\in\mathbb{N}_{0}} \bigl[ K_{n,\widehat{\pi}%
_{n-1}^{ ( m ) }}^{i} ( f_{n} ) -\omega_{n} \bigl( \widehat{%
\pi}_{n-1}^{ ( m ) } \bigr) ( f_{n} ) \bigr] \biggr)
\biggr\vert\nonumber\\
&&\qquad =\sum_{i\in\mathbb{N}_{0}} \bigl\vert\bigl( \omega_{n} \bigl( \widehat{\pi}%
_{n-1}^{ ( m+1 ) } \bigr) -\omega_{n} \bigl( \widehat{\pi}%
_{n-1}^{ ( m ) } \bigr) \bigr) K_{n,\widehat{\pi}_{n-1}^{ (
m ) }}^{i} ( f_{n} ) \bigr\vert\\
&&\qquad\leq\sum_{i\in\mathbb{N}_{0}}\rho_{n}^{i} \bigl\Vert\omega_{n} \bigl(
\widehat{\pi}_{n-1}^{ ( m+1 ) } \bigr) -\omega_{n} \bigl( \widehat{%
\pi}_{n-1}^{ ( m ) } \bigr) \bigr\Vert_{\mathrm{tv}} \nonumber
\nonumber\\
&&\qquad \leq\frac{2}{ ( 1-\rho_{n} ) ^{2}} \bigl\Vert\widehat{\pi}%
_{n-1}^{ ( m+1 ) }-\widehat{\pi}_{n-1}^{ ( m )
} \bigr\Vert_{\mathrm{tv}},\nonumber
\end{eqnarray}
using Assumption~\ref{assumA1} and (\ref
{eq:sensitivityinvariantdistribution}) in
Proposition~\ref{proposition:sensitivity}.

Now for any $f_{n-1}\in B ( E_{n-1} ) $, we have%
\[
\widehat{\pi}_{n-1}^{ ( m+1 ) } ( f_{n-1} ) -\widehat{\pi}%
_{n-1}^{ ( m ) } ( f_{n-1} ) =\frac{f_{n-1} ( \mathbf{X%
}_{n-1}^{ ( m+1 ) } ) }{m+2}-\frac{\widehat{\pi}_{n-1}^{ (
m ) } ( f_{n-1} ) }{m+2},
\]
thus
%
%e6.15 ###
%
\begin{equation}\label{eq:TVempiricalmeasures}
\bigl\Vert\widehat{\pi}_{n-1}^{ ( m+1 ) }-\widehat{\pi}%
_{n-1}^{ ( m ) } \bigr\Vert_{\mathrm{tv}}\leq\frac{2}{m+2}.
\end{equation}
The result follows now directly combining (\ref%
{eq:identitypoissonequationdelmoralmiclo}), (\ref
{eq:boundonfirstterm}), (%
\ref{eq:boundon2ndterm}), (\ref{eq:TVempiricalmeasures}) and using
Minkowski's inequality.
\end{pf}
\begin{proposition}
\label{proposition:firstpartergodicaverages}Given Assumption \ref
{assumA1}, for any $n\in\{
2,\ldots,P \} $ and any $p\geq1$ there exists $B_{1,n},B_{2,p}<\infty$
such that for $\mathbf{x}_{1\dvtx n}^{ ( 0 ) },f_{n}\in B (
E_{n} ) $ and $i\in\mathbb{N}_{0}$%
\[
\mathbb{E}_{\mathbf{x}_{1\dvtx n}^{ ( 0 ) }} \bigl[ \bigl\vert\widehat{%
\pi}_{n}^{ ( i ) } ( f_{n} ) -S_{n}^{ ( i )
} ( f_{n} ) \bigr\vert^{p} \bigr] ^{1/p}\leq\frac{B_{1,n}B_{2,p}%
}{ ( i+1 ) ^{{1/2}}}.
\]
\end{proposition}
\begin{pf}
Using (\ref{eq:longdecompositionaverages}) and Minkowski's
inequality, we obtain%
%
%e6.16 ###
%
\begin{eqnarray}\label{eq:inequalitydifference}\qquad
&& \mathbb{E}_{\mathbf{x}_{1\dvtx n}^{ ( 0 ) }} \bigl[ \bigl\vert\widehat{%
\pi}_{n}^{ ( i ) } ( f_{n} ) -S_{n}^{ ( i )
} ( f_{n} ) \bigr\vert^{p} \bigr] ^{1/p}
\nonumber\\[-2pt]
&&\qquad \leq\frac{1}{ ( i+1 ) }\mathbb{E}_{\mathbf{x}_{1\dvtx n}^{ (
0 ) }} \bigl[ \bigl\vert M_{n}^{ ( i+1 ) } ( f_{n} )
\bigr\vert^{p} \bigr] ^{1/p} \nonumber\\[-2pt]
&&\qquad\quad{} +\frac{1}{ ( i+1 ) }\sum_{m=0}^{i}\mathbb{E}_{\mathbf{x}%
_{1\dvtx n}^{ ( 0 ) }} \bigl[ \bigl\vert\widehat{f}_{n,\widehat{\pi}%
_{n-1}^{ ( m+1 ) }} \bigl( \mathbf{X}_{n}^{ ( m+1 )
} \bigr) -\widehat{f}_{n,\widehat{\pi}_{n-1}^{ ( m ) }} \bigl(
\mathbf{X}_{n}^{ ( m+1 ) } \bigr) \bigr\vert^{p} \bigr] ^{1/p}
\\[-2pt]
&&\qquad\quad{} +\frac{1}{i+1}\mathbb{E}_{\mathbf{x}_{1\dvtx n}^{ ( 0 ) }} \bigl[
\bigl\vert\widehat{f}_{n,\widehat{\pi}_{n-1}^{ ( 0 ) }} \bigl(
\mathbf{X}_{n}^{ ( 0 ) } \bigr) \bigr\vert^{{p}} \bigr] ^{1/p}\nonumber\\[-2pt]
&&\qquad\quad{}+
\frac{1}{i+1}\mathbb{E}_{\mathbf{x}_{1\dvtx n}^{ ( 0 ) }} \bigl[
\bigl\vert\widehat{f}_{n,\widehat{\pi}_{n-1}^{ ( i+1 ) }} \bigl(
\mathbf{X}_{n}^{ ( i+1 ) } \bigr) \bigr\vert^{{p}} \bigr] ^{1/p}.
\nonumber
\end{eqnarray}
The first term on the right-hand side of (\ref{eq:inequalitydifference}) is bounded
using Proposition~\ref{proposition:martingalebounded}, the term on the last
line of the right-hand side are going to zero because of
Proposition~\ref{proposition:poissonbounded}. For the second term, we obtain using
Proposition~\ref{proposition:poissonincrements}
\begin{eqnarray*}
\sum_{m=0}^{i}\mathbb{E}_{\mathbf{x}_{1\dvtx n}^{ ( 0 ) }} \bigl[
\bigl\vert\widehat{f}_{n,\widehat{\pi}_{n-1}^{ ( m+1 ) }} \bigl(
\mathbf{X}_{n}^{ ( m+1 ) } \bigr) -\widehat{f}_{n,\widehat{\pi}%
_{n-1}^{ ( m ) }} \bigl( \mathbf{X}_{n}^{ ( m+1 ) } \bigr)
\bigr\vert^{p} \bigr] ^{1/p} & \leq & \sum_{m=0}^{i}\frac{B_{n}}{m+2} \\[-2pt]
& \leq & B_{n}\log( i+2 ).
\end{eqnarray*}
The result follows.
\end{pf}
\begin{pf*}{Proof of Theorem~\ref{theorem:convergenceofaverages}}
Under Assumption~\ref{assumA1}, the
result is clearly true for $n=1$ thanks to Lemma \ref%
{lemma:k1geometricallyergodic}. Assume it is true for $n-1$ and let us prove
it for $n$. We have, using Minkowski's inequality,
\begin{eqnarray*}
\mathbb{E}_{\mathbf{x}_{1\dvtx n}^{ ( 0 ) }} \bigl[ \bigl\vert\widehat{%
\pi}_{n}^{ ( i ) } ( f_{n} ) -\pi_{n} ( f_{n} )
\bigr\vert^{p} \bigr] ^{1/p} &\leq& \mathbb{E}_{\mathbf{x}_{1\dvtx n}^{ (
0 ) }} \bigl[ \bigl\vert\widehat{\pi}_{n}^{ ( i ) } (
f_{n} ) -S_{n}^{ ( i ) } ( f_{n} ) \bigr\vert^{p}%
\bigr] ^{1/p} \\[-2pt]
&&{} + \mathbb{E}_{\mathbf{x}_{1\dvtx n}^{ ( 0 ) }} \bigl[ \bigl\vert
S_{n}^{ ( i ) } ( f_{n} ) -\pi_{n} ( f_{n} )
\bigr\vert^{p} \bigr] ^{1/p}.
\end{eqnarray*}
The first term on the right-hand side can be bounded using Proposition
\ref{proposition:firstpartergodicaverages}. For the second term, we have
\begin{eqnarray*}
&&\mathbb{E}_{\mathbf{x}_{1\dvtx n}^{ ( 0 ) }} \bigl[ \bigl\vert
S_{n}^{ ( i ) } ( f_{n} ) -\pi_{n} ( f_{n} )
\bigr\vert^{p} \bigr] ^{1/p}\\[-2pt]
&&\qquad\leq\frac{1}{ ( i+1 ) }\sum_{m=0}^{i}%
\mathbb{E}_{\mathbf{x}_{1\dvtx n}^{ ( 0 ) }} \bigl[ \bigl\vert\omega
_{n} \bigl( \widehat{\pi}_{n-1}^{ ( m ) } \bigr) (
f_{n} ) -\omega_{n} ( \pi_{n-1} ) ( f_{n} )
\bigr\vert^{p} \bigr] ^{1/p}.
\end{eqnarray*}
Using (\ref{eq:invariantdistributionperturbedkernel}), we obtain%
\begin{eqnarray*}
&& \omega_{n} ( \pi_{n-1} ) ( f_{n} ) -\omega_{n} \bigl(
\widehat{\pi}_{n-1}^{ ( m ) } \bigr) ( f_{n} ) \\[-2pt]
&&\qquad =\frac{ ( \pi_{n-1}\times\overline{\pi}_{n} ) ( \pi
_{n/n-1}\cdot f_{n} ) }{\pi_{n-1} ( \pi_{n/n-1} ) }-\frac{ (
\widehat{\pi}_{n-1}^{ ( m ) }\times\overline{\pi}_{n} )
( \pi_{n/n-1}\cdot f_{n} ) }{\widehat{\pi}_{n-1}^{ ( m )
} ( \pi_{n/n-1} ) } \\
&&\qquad =\frac{ ( ( \pi_{n-1}-\widehat{\pi}_{n-1}^{ ( m )
} ) \times\overline{\pi}_{n} ) ( \pi_{n/n-1}\cdot f_{n} ) \cdot %
\widehat{\pi}_{n-1}^{ ( m ) } ( \pi_{n/n-1} ) }{\widehat{%
\pi}_{n-1}^{ ( m ) } ( \pi_{n/n-1} ) \cdot
\pi_{n-1} ( \pi_{n/n-1} ) } \\
&&\qquad\quad{} +\frac{ ( \widehat{\pi}_{n-1}^{ ( m ) }\times\overline{\pi}%
_{n} ) ( \pi_{n/n-1}\cdot f_{n} ) \cdot  ( \widehat{\pi}%
_{n-1}^{ ( m ) }-\pi_{n-1} ) ( \pi_{n/n-1} ) }{%
\widehat{\pi}_{n-1}^{ ( m ) } ( \pi_{n/n-1} ) \cdot
\pi_{n-1} ( \pi_{n/n-1} ) }
\end{eqnarray*}
so, as $\pi_{n-1} ( \pi_{n/n-1} ) =1$, we have
\begin{eqnarray*}
&& \bigl\vert\omega_{n} ( \pi_{n-1} ) ( f_{n} ) -\omega
_{n} \bigl( \widehat{\pi}_{n-1}^{ ( m ) } \bigr) (
f_{n} ) \bigr\vert\\
&&\qquad \leq\bigl\vert\bigl( \bigl( \pi_{n-1}-\widehat{\pi}_{n-1}^{ (
m ) } \bigr) \times\overline{\pi}_{n} \bigr) ( \pi
_{n/n-1}\cdot f_{n} ) \bigr\vert\\
&&\qquad\quad{} +\frac{ \vert( \widehat{\pi}_{n-1}^{ ( m ) }\times
\overline{\pi}_{n} ) ( \pi_{n/n-1}\cdot f_{n} ) \cdot  ( \widehat{%
\pi}_{n-1}^{ ( m ) }-\pi_{n-1} ) ( \pi_{n/n-1} )
\vert}{\widehat{\pi}_{n-1}^{ ( m ) } ( \pi
_{n/n-1} ) }.
\end{eqnarray*}
Assumption~\ref{assumA1} implies that there exists $D_{n}<\infty$ such
that $%
\pi_{n/n-1} ( \mathbf{x}_{n-1} ) <D_{n}$ over $E_{n-1}$. Thus, we
have using the induction hypothesis%
\begin{eqnarray*}
&& \mathbb{E}_{\mathbf{x}_{1\dvtx n}^{ ( 0 ) }} \bigl[\bigl \vert\omega
_{n} \bigl( \widehat{\pi}_{n-1}^{ ( m ) } \bigr) (
f_{n} ) -\omega_{n} ( \pi_{n-1} ) ( f_{n} )
\bigr\vert^{p} \bigr] ^{1/p} \\
&&\qquad \leq2D_{n}\mathbb{E}_{\mathbf{x}_{1\dvtx n}^{ ( 0 ) }} \biggl[
\biggl\vert\widehat{\pi}_{n-1}^{ ( m ) } \biggl( \frac{\pi_{n/n-1}}{%
D_{n}} \biggr) -\pi_{n-1} \biggl( \frac{\pi_{n/n-1}}{D_{n}} \biggr)
\biggr\vert^{p} \biggr] ^{1/p} \\
&&\qquad \leq\frac{2D_{n}C_{1,n-1}C_{2,p}}{ ( m+1 ) ^{1/2}}
\end{eqnarray*}
and
\begin{eqnarray*}
\mathbb{E}_{\mathbf{x}_{1\dvtx n}^{ ( 0 ) }} \bigl[ \bigl\vert
S_{n}^{ ( i ) } ( f_{n} ) -\pi_{n} ( f_{n} )
\bigr\vert^{p} \bigr] ^{1/p} &\leq&
\frac{2D_{n}C_{1,n-1}C_{2,p}}{ (
i+1 ) }\sum_{m=0}^{i}\frac{1}{ ( m+1 ) ^{1/2}}\\
&\leq&\frac{D_{n}C_{1,n-1}C_{2,p}}{ ( i+1 ) ^{1/2}}.
\end{eqnarray*}
This concludes the proof.
\end{pf*}

%s6.4 ###
\subsection{Convergence of marginals}

\mbox{}

\begin{pf*}{Proof of Theorem~\ref{theorem:convergenceofmarginals}}
The proof is adapted from Proposition 4 in~\cite{andrieumoulines2006}. For $n=1$ the result follows directly from
Assumption
\ref{assumA1}. Now consider the case where $n\geq2$. We use the following
decomposition for $0\leq n ( i ) \leq i$:
\begin{eqnarray*}
&&
\bigl\vert\mathbb{E}_{\mathbf{x}_{1\dvtx n}^{ ( 0 ) }} \bigl[
f_{n} \bigl( \mathbf{X}_{n}^{ ( i ) } \bigr) -\pi_{n} (
f_{n} ) \bigr] \bigr\vert\\
&&\qquad \leq\bigl\vert\mathbb{E}_{\mathbf{x}%
_{1\dvtx n}^{ ( 0 ) }} \bigl[ f_{n} \bigl( \mathbf{X}_{n}^{ ( i )
} \bigr) -K_{n,\widehat{\pi}_{n-1}^{ ( i-n ( i ) )
}}^{n ( i ) }f_{n} \bigl( \mathbf{X}_{n}^{ ( i-n ( i )
) } \bigr) \bigr] \bigr\vert\\
&&\qquad\quad{} + \bigl\vert\mathbb{E}_{\mathbf{x}_{1\dvtx n}^{ ( 0 ) }} \bigl[ K_{n,%
\widehat{\pi}_{n-1}^{ ( i-n ( i ) ) }}^{n ( i )
} ( f_{n} ) \bigl( \mathbf{X}_{n}^{ ( i-n ( i )
) } \bigr) -\omega_{n} \bigl( \widehat{\pi}_{n-1}^{ ( i-n (
i ) ) } \bigr) ( f_{n} ) \bigr] \bigr\vert\\
&&\qquad\quad{} + \bigl\vert\mathbb{E}_{\mathbf{x}_{1\dvtx n}^{ ( 0 ) }} \bigl[ \omega
_{n} \bigl( \widehat{\pi}_{n-1}^{ ( i-n ( i ) ) } \bigr)
( f_{n} ) -\omega_{n} ( \pi_{n-1} ) ( f_{n} ) %
\bigr] \bigr\vert.
\end{eqnarray*}
Assumption~\ref{assumA1} implies that
\[
\bigl\vert\mathbb{E}_{\mathbf{x}_{1\dvtx n}^{ ( 0 ) }} \bigl[ K_{n,%
\widehat{\pi}_{n-1}^{ ( i-n ( i ) ) }}^{n ( i )
} ( f_{n} ) \bigl( \mathbf{X}_{n}^{ ( i-n ( i )
) } \bigr) -\omega_{n} \bigl( \widehat{\pi}_{n-1}^{ ( i-n (
i ) ) } \bigr) ( f_{n} ) \bigr] \bigr\vert\leq\rho
_{n}^{n ( i ) }.
\]
For the first term, we use the following decomposition:
\begin{eqnarray*}
&& \mathbb{E}_{\mathbf{x}_{1\dvtx n}^{ ( 0 ) }} \bigl[ f_{n} \bigl( \mathbf{%
X}_{n}^{ ( i ) } \bigr) -K_{n,\widehat{\pi}_{n-1}^{ (
i-n ( i ) ) }}^{n ( i ) } ( f_{n} ) \bigl(
\mathbf{X}_{n}^{ ( i-n ( i ) ) } \bigr) \bigr] \\
&&\qquad =\sum_{j=2}^{n ( i ) }\mathbb{E}_{\mathbf{x}_{1\dvtx n}^{ (
0 ) }} \bigl[ K_{n,\widehat{\pi}_{n-1}^{ ( i-j+1 )
}}^{j-1} ( f_{n} ) \bigl( \mathbf{X}_{n}^{ ( i-j+1 )
} \bigr) -K_{n,\widehat{\pi}_{n-1}^{ ( i-j ) }}^{j} (
f_{n} ) \bigl( \mathbf{X}_{n}^{ ( i-j ) } \bigr) \bigr]
\end{eqnarray*}
and
\begin{eqnarray*}
&& \mathbb{E}_{\mathbf{x}_{1\dvtx n}^{ ( 0 ) }} \bigl[ K_{n,\widehat{\pi}%
_{n-1}^{ ( i-j+1 ) }}^{j-1} ( f_{n} ) \bigl( \mathbf{X}%
_{n}^{ ( i-j+1 ) } \bigr) -K_{n,\widehat{\pi}_{n-1}^{ (
i-j ) }}^{j} ( f_{n} ) \bigl( \mathbf{X}_{n}^{ (
i-j ) } \bigr) \bigr] \\
&&\qquad =\mathbb{E}_{\mathbf{x}_{1\dvtx n}^{ ( 0 ) }} \bigl[ \mathbb{E}_{%
\mathbf{x}_{1\dvtx n}^{ ( 0 ) }} \bigl[ K_{n,\widehat{\pi}%
_{n-1}^{ ( i-j+1 ) }}^{j-1} ( f_{n} ) \bigl( \mathbf{X}%
_{n}^{ ( i-j+1 ) } \bigr)\\
&&\qquad\quad\hspace*{53.5pt}{} -K_{n,\widehat{\pi}_{n-1}^{ (
i-j ) }}^{j-1} ( f_{n} ) \bigl( \mathbf{X}_{n}^{ (
i-j+1 ) } \bigr) \vert\mathcal{G}_{n}^{i-j} \bigr] \bigr],
\end{eqnarray*}
where
\begin{eqnarray*}
&& K_{n,\widehat{\pi}_{n-1}^{ ( i-j+1 ) }}^{j-1} ( f_{n} )
\bigl( \mathbf{X}_{n}^{ ( i-j+1 ) } \bigr) -K_{n,\widehat{\pi}%
_{n-1}^{ ( i-j ) }}^{j-1} ( f_{n} ) \bigl( \mathbf{X}%
_{n}^{ ( i-j+1 ) } \bigr) \\
&&\qquad =\sum_{m=0}^{j-2}K_{n,\widehat{\pi}_{n-1}^{ ( i-j+1 )
}}^{m} ( K_{n,\widehat{\pi}_{n-1}^{ ( i-j+1 ) }}-K_{n,\widehat{%
\pi}_{n-1}^{ ( i-j ) }} ) K_{n,\widehat{\pi}_{n-1}^{ (
i-j ) }}^{j-1-m-1} ( f_{n} ) \bigl( \mathbf{X}_{n}^{ (
i-j+1 ) } \bigr) \\
&&\qquad =\sum_{m=0}^{j-2}K_{n,\widehat{\pi}_{n-1}^{ ( i-j+1 )
}}^{m} ( K_{n,\widehat{\pi}_{n-1}^{ ( i-j+1 ) }}-K_{n,\widehat{%
\pi}_{n-1}^{ ( i-j ) }} ) \\
&&\qquad\quad\hspace*{15.6pt}{} \times\bigl( K_{n,\widehat{\pi}_{n-1}^{ ( i-j )
}}^{j-1-m-1} ( f_{n} ) \bigl( \mathbf{X}_{n}^{ ( i-j+1 )
} \bigr) -\omega_{n} \bigl( \widehat{\pi}_{n-1}^{ ( i-j ) } \bigr)
( f_{n} ) \bigr) .
\end{eqnarray*}
Now we have from Proposition~\ref{proposition:sensitivity} that
\begin{eqnarray*}
\Vert K_{n,\widehat{\pi}_{n-1}^{ ( i-j+1 ) }}^{m} (
\mathbf{x}_{n},\cdot) -K_{n,\widehat{\pi}_{n-1}^{ ( i-j )
}}^{m} ( \mathbf{x}_{n},\cdot) \Vert_{\mathrm{tv}}&
\leq&\frac{2}{ ( 1-\rho_{n} ) } \bigl\Vert\widehat{\pi}%
_{n-1}^{ ( i-j+1 ) }-\widehat{\pi}_{n-1}^{ ( i-j )
} \bigr\Vert_{\mathrm{tv}} \\
& \leq&\frac{2}{ ( 1-\rho_{n} ) }\frac{1}{i-j+2}
\end{eqnarray*}
and using Assumption~\ref{assumA1}
\begin{eqnarray*}
&& \bigl\vert\mathbb{E}_{\mathbf{x}_{1\dvtx n}^{ ( 0 ) }} \bigl[ \mathbb{E%
}_{\mathbf{x}_{1\dvtx n}^{ ( 0 ) }} \bigl[ K_{n,\widehat{\pi}%
_{n-1}^{ ( i-j+1 ) }}^{j-1} ( f_{n} ) \bigl( \mathbf{X}%
_{n}^{ ( i-j+1 ) } \bigr) -K_{n,\widehat{\pi}_{n-1}^{ (
i-j ) }}^{j-1} ( f_{n} ) \bigl( \mathbf{X}_{n}^{ (
i-j+1 ) } \bigr) \vert\mathcal{G}_{n}^{i-j} \bigr] \bigr]
\bigr\vert\\
&&\qquad \leq\Biggl\vert\mathbb{E}_{\mathbf{x}_{1\dvtx n}^{ ( 0 ) }} \Biggl[
\mathbb{E}_{\mathbf{x}_{1\dvtx n}^{ ( 0 ) }} \Biggl[ \sum_{m=0}^{j-2}K_{n,%
\widehat{\pi}_{n-1}^{ ( i-j+1 ) }}^{m} ( K_{n,\widehat{\pi}%
_{n-1}^{ ( i-j+1 ) }}-K_{n,\widehat{\pi}_{n-1}^{ ( i-j )
}} ) \\
&&\qquad\quad\hspace*{78.7pt}{} \times\bigl( K_{n,\widehat{\pi}%
_{n-1}^{ ( i-j ) }}^{j-1-m-1} ( f_{n} ) \bigl( \mathbf{X}%
_{n}^{ ( i-j+1 ) } \bigr) \\
&&\qquad\quad\hspace*{130.2pt}{}-\omega_{n} \bigl( \widehat{\pi}%
_{n-1}^{ ( i-j ) } \bigr) ( f_{n} ) \bigr) \vert
\mathcal{G}_{n}^{i-j} \Biggr] \Biggr] \Biggr\vert\\
&&\qquad \leq\frac{2}{ ( 1-\rho_{n} ) ( i-j+2 ) }%
\sum_{m=0}^{j-2}\rho_{n}^{j-m-2}
=\frac{2}{ ( 1-\rho_{n} ) ( i-j+2 ) }\frac{1-\rho
_{n}^{j-1}}{1-\rho_{n}}
\end{eqnarray*}
and
\begin{eqnarray*}
&& \bigl\vert\mathbb{E}_{\mathbf{x}_{1\dvtx n}^{ ( 0 ) }} \bigl[
f_{n} \bigl( \mathbf{X}_{n}^{ ( i ) } \bigr) -K_{n,\widehat{\pi}%
_{n-1}^{ ( i-n ( i ) ) }}^{n ( i ) } (
f_{n} ) \bigl( \mathbf{X}_{n}^{ ( i-n ( i ) )
} \bigr) \bigr] \bigr\vert\\
&&\qquad \leq\frac{2}{ ( 1-\rho_{n} ) ^{2}}\sum_{j=2}^{n ( i ) }%
\frac{1}{ ( i-j+2 ) } \\
&&\qquad \leq\frac{2}{ ( 1-\rho_{n} ) ^{2}}\log\biggl( \frac{i}{%
i-n ( i ) +1} \biggr) .
\end{eqnarray*}
Finally, to study the last term $\mathbb{E} [ \omega_{n} ( \widehat{%
\pi}_{n-1}^{ ( i-n ( i ) ) } ) ( f_{n} )
-\omega_{n} ( \pi_{n-1} ) ( f_{n} ) ] $, we use
the same decomposition used in the proof of Theorem \ref%
{theorem:convergenceofaverages} to obtain
\begin{eqnarray*}
&& \bigl\vert\mathbb{E} \bigl[ \omega_{n} \bigl( \widehat{\pi}_{n-1}^{ (
i-n ( i ) ) } \bigr) ( f_{n} ) -\omega_{n} (
\pi_{n-1} ) ( f_{n} ) \bigr] \bigr\vert\\
&&\qquad \leq2D_{n}\mathbb{E}_{\mathbf{x}_{1\dvtx n}^{ ( 0 ) }} \biggl[
\biggl\vert\widehat{\pi}_{n-1}^{ ( i-n ( i ) ) } \biggl(
\frac{\pi_{n/n-1}}{D_{n}} \biggr) -\pi_{n-1} \biggl( \frac{\pi_{n/n-1}}{D_{n}%
} \biggr) \biggr\vert\biggr] \\
&&\qquad \leq\frac{2D_{n}C_{1,n-1}}{ ( i-n ( i ) +1 ) ^{1/2}}.
\end{eqnarray*}
One can\vspace*{1pt} check that $ \vert\mathbb{E}_{\mathbf{x}_{1\dvtx n}^{ ( 0 )
}} [ f_{n} ( \mathbf{X}_{n}^{ ( i ) } ) -\pi
_{n} ( f_{n} ) ] \vert$ converges toward zero for $%
n ( i ) = \lfloor i^{\alpha} \rfloor$ where $0<\alpha
<1$.
\end{pf*}
\end{appendix}

%suskaldyti doi

%
\printaddresses

\end{document}